%% file: main.tex
\documentclass[letterpaper, paper,11pt]{AAS}		% for final proceedings (20-page limit)

\input{preamble}

\title{Leveraging Hamiltonian Structure for Accurate Uncertainty Propagation}
\author{Amit Jain\footnote{Ph.D. Candidate, Department of Aerospace Engineering, Pennsylvania State University, State College, PA-16802, Email: axj307@psu.edu.}, Puneet Singla\footnote{Professor, AIAA Associate Fellow, AAS Fellow, Department of Aerospace Engineering, Pennsylvania State University, State College, PA-16802, Email: psingla@psu.edu.}, and Roshan Eapen\footnote{Assistant Professor, AIAA Member, AAS Member, Department of Aerospace Engineering, Pennsylvania State University, State College, PA-16802, Email: reapen@psu.edu.}}

\PaperNumber{23-361}

\begin{document}

\maketitle

\begin{abstract} 
In this work, we leverage the Hamiltonian kind structure for accurate uncertainty propagation through a nonlinear dynamical system. The developed approach utilizes the fact that the stationary probability density function is purely a function of the Hamiltonian of the system. This fact is exploited to define the basis functions for approximating the solution of the Fokker-Planck-Kolmogorov equation. This approach helps in curtailing the growth of basis functions with the state dimension. Furthermore, sparse approximation tools have been utilized to automatically select appropriate basis functions from an over-complete dictionary. A nonlinear oscillator and two-body problem are considered to show the efficacy of the proposed approach. Simulation results show that such an approach is effective in accurately propagating uncertainty through non-conservative as well as conservative systems.

% A computationally efficient approach is presented to compute the reachability set for a nonlinear system. The central idea of the developed approach is to represent the reachability set as the probability density function (pdf) and find the evolution of the state pdf through a governing equation known as Fokker-Plank-Kolmogorov Equation (FPKE). A non-product sampling method known as Conjugate Unscented Transformation (CUT) is used as collocation points to alleviate the curse of dimensionality. In conjunction with the minimal cubature points provided by the CUT method, an $l_1$ norm minimization technique is utilized to optimally select the appropriate basis functions from a larger complete dictionary of basis functions. The main contribution of this paper is to include Hamiltonian basis functions in an over-complete dictionary of monomial basis functions and utilize a sparse algorithm to analyze the dependency of the stationary pdf on Hamiltonians. The monomial basis functions are necessary to approximate the transient pdf, while the Hamiltonian basis functions dictate the stationary behavior. Various numerical examples are shown to demonstrate the utility of the proposed approach.

\end{abstract}

\section{Introduction}

\input{src/01_intro}

\section{Problem Formulation}

\input{src/02_prob}

\section{Numerical Simulations}
\input{src/03_numeric}

\section{Conclusion}

\input{src/04_conclusion}

\section*{Acknowledgments}
This material is based upon work supported through United States Air Force Office of Scientific Research (AFOSR) grant FA9550-20-1-0176.

%Additional information:
%PSD plot filters, 
% \newpage
\bibliographystyle{AAS_publication.bst}   % Number the references.
\bibliography{reach.bib,allrefs.bib}

\end{document}

%% file: preamble.tex
\usepackage{new-mysymb}
\usepackage{longtable}
\usepackage{multirow}
\usepackage{pbox}
\usepackage[hyphens]{url}
\usepackage{soul}
\usepackage{nomencl}
\usepackage{etoolbox}

\usepackage{rotating}
\usepackage[BF,bf]{subfigure}
\usepackage{mathtools}

\usepackage{eqlist} % Makes for a nice list of symbols.
\usepackage[nosepfour,warning,np,debug,autolanguage]{numprint}
\usepackage{acro}
\usepackage{booktabs}
\usepackage{graphicx}
\usepackage{figsize}
\usepackage{float}
\usepackage{textcomp}
\usepackage{url}
\usepackage{colortbl}

\definecolor{orange}{cmyk}{0,0.3,1,0}

\renewcommand\nomgroup[1]{%
   \item[\bfseries
   \ifstrequal{#1}{A}{Subscripts / Superscripts}]}

\makenomenclature
\printnomenclature

\usepackage[utf8]{inputenc}
\usepackage{amsmath}
\usepackage[version=4]{mhchem}
\usepackage{siunitx}
\usepackage{tabularx}
\usepackage{caption}
\usepackage{listings}
\usepackage{graphicx,color,subfigure,wrapfig,times,breqn}
\usepackage[colorlinks=true, pdfstartview=FitV, linkcolor=black, citecolor= black, urlcolor= black]{hyperref}
\usepackage[section]{placeins}
\usepackage{url}
\usepackage{enumitem}
\usepackage{subfigure}

\usepackage{algorithm,algcompatible,algorithmicx}
\algnewcommand\INPUT{\item[\textbf{Input:}]}%
\algnewcommand\OUTPUT{\item[\textbf{Output:}]}%

\usepackage{array}

\usepackage{subcaption,mwe}

\usepackage{array}
\newcolumntype{M}[1]{>{\centering\arraybackslash}m{#1}}

%% file: src/01_intro.tex
Analysis of an engineering system relies heavily on the quantification and propagation of uncertainty. Uncertainty Propagation (UP) of the state pdf through a nonlinear dynamical model is a complex process. A simplistic approach for uncertainty propagation, Monte Carlo (MC) method, involves evaluating the governing equations with many different realizations of the uncertain parameter and appropriately deriving the statistics from the outputs. MC methods require extensive computational resources and effort and become increasingly infeasible for high-dimensional dynamic systems \cite{herrador2005estimation}. Furthermore, MC methods provide only an approximate description of the uncertainty propagation problem by restricting the solution to a small number of parameters - such as the first N moments of the sought state pdf. 

% The response of a nonlinear dynamical system to a Gaussian white noise excitation is a diffusion process and the transition pdf of the response is governed by the Fokker-Planck-Kolmogrov equation (FPKE). Although FPKE is a linear partial differential equation, it is a formidable equation to solve. The only way to obtain the exact solutions is by solving the FPKE under appropriate initial and boundary conditions. Closed-form solutions of FPKE exist only for a few cases, mainly demonstrated in a linear system under the influence of additive white Gaussian noise \cite{}. The exact stationary solution for few classes of nonlinear stochastic systems has been shown by Fuller \cite{}. These stationary solutions are closely related to the Maxwell-Boltzmann distributon in classical statistical mechanics. For a stochastically excited and dissipated Hamiltonian system, this property is manifested by that the stationary pdf which is a function of the Hamiltonian.  

The time evolution of the state pdf for a nonlinear dynamical system driven by the Gaussian white noise process is governed by the Fokker-Planck-Kolmogorov Equation (FPKE) \cite{Fuller1969}. The FPKE is a convection-diffusion partial differential equation, where convection corresponds to the drift of the mean response and diffusion corresponds to the dispersion or scattering of the family of responses. Even though FPKE is a linear partial differential equation, it is a complex problem to solve for general nonlinear systems. In order to obtain an exact solution, the FPKE must be solved under appropriate initial and boundary conditions. Closed-form solutions of FPKE exist only for a few cases, mainly demonstrated in a linear system under the influence of additive white Gaussian noise. The solution of FPKE for linear systems in conjunction with Bayes’s rule to assimilate measurement data with a linear model leads to the celebrated Kalman Filter \cite{kumar2006multi}. The exact stationary solution for a few classes of nonlinear stochastic systems has been shown in \cite{Fuller1969,caughey1982exact,cai1996exact}. These stationary solutions are closely related to the Maxwell-Boltzmann distribution in classical statistical mechanics.

Although there has been some success in finding exact solutions, it has been necessary to adopt approximate numerical solution procedures for a general nonlinear system. Numerous techniques such as the variational methods \cite{jordan1998variational,johnson1980extension}, the finite element method \cite{spencer1993numerical}, the finite difference method \cite{roberts1986first,cichon2003approximation}, path integration method \cite{wehner1983numerical,xu2019path,naess1992response}, and weighted residuals schemes \cite{soize1988steady,wen1975approximate} have all been developed as approximate numerical techniques. Methods such as weighted residual schemes involve the computation of high-dimensional integrals and, thus, suffer from the curse of dimensionality. Finite Element and finite difference-based methods are severely handicapped for higher dimensions because the discretization of the space over which pdf lives is computationally impractical. The variational method uses an eigenfunction expansion of the FPKE and involves satisfying the complex boundary conditions, which are hard to satisfy. The path integral or cell mapping method is highly efficient for 2-D systems but is computationally very expensive for higher dimension systems. Recently, much attention has been directed towards Adaptive Gaussian Mixture Model (AGMM) methods to solve the FPKE \cite{terejanu2008uncertainty} accurately. The key idea of the AGMM is to approximate the state pdf by a finite sum of Gaussian density functions, and the accuracy of these methods relies on optimally selecting the number and location of mixture model components. Furthermore, the solution to the convex optimization problem involves the computation of multi-dimensional integrals, which can be a computationally expensive affair for high-dimensional systems.

In summary, the key limitation in computing the solution to the FPKE lies in the curse of dimensionality, which stems from the discretization of the domain of the pdf and/or the computation of multi-dimensional projection integrals. The principal contribution of this work is to address this computation challenge. Recently developed sparse collocation methods \cite{mercurio2016conjugate,mirzaei2018sparse} have been used to solve the FPKE efficiently. The solution process involves the finite series expansion of the log-pdf in terms of suitable polynomial basis functions. The coefficient and order of the finite series expansion for log-pdf are determined by exactly satisfying the FPKE at the collocation points. 

The main challenge in the development of any collocation method lies in choosing the appropriate collocation points and the basis functions. The first challenge is taken care of by selecting the minimal number of collocation points required to approximate the domain accurately. Thus, a non-product sampling method known as Conjugate Unscented Transformation (CUT) is used as collocation points. The second challenge is the selection of the appropriate basis function dictionary, which is the major contribution of this paper. Exploiting the Hamiltonian structure, it can be found that the log-pdf of the stationary solution is directly proportional to the Hamiltonian function. This result can be directly obtained by solving the FPKE. Thus, to capture the pdf, the basis function dictionary must also include Hamiltonians in addition to the commonly used monomials. By assuming an over-complete dictionary of monomials and Hamiltonians basis function, the sparse approximation technique can nicely capture the Hamiltonian coefficients as the system reaches the stationary value. Therefore, the main contribution of this paper is including the Hamiltonians in the basis functions dictionary to compute the uncertainty propagation. 

% A novel idea to exploit the Hamiltonian structure leads to the effective
% The major contribution of this paper is the selection of the basis function dictionary in approximating the log-pdf at each time instance. 
% The second challenge of selecting the appropriate basis function dictionary is taken care of by exploiting the Hamiltonian function. 
% can be combat to some extent by selecting the basis function dictionary in an effective manner. 

The structure of the paper is as follows: first, a brief overview of the FPKE is presented, followed by motivation and detailed development of the proposed methodology. The constraints of the solution process are discussed in detail, including the selection of the collocation points and the determination of a minimal order interpolation polynomial for the log-pdf. Numerical examples, including a nonlinear oscillator and a two-body problem, are then presented to demonstrate the utility of the developed approach.

%% file: src/02_prob.tex
% This section computes the reachability set for dynamical systems represented by the \eqref{Eq:StochasticSystem}. A reachability envelope or set is defined as the computation of the target state pdf generated due to the uncertainties present in the system.
This section computes the uncertainty propagation for dynamical systems represented by the \eqref{Eq:StochasticSystem}. These uncertainties can be introduced due to errors in the initial condition, external control input, or modeling errors in the system. Consider the following stochastic dynamical system:
\begin{equation}\label{Eq:StochasticSystem}
d \mathbf{x}(t) = \mathbf{f}(\mathbf{x}(t), t) dt + \mathbf{G}(\mathbf{x}(t), t) d \Gamma(t) 
\end{equation}
where $\x \in \mathbb{R}^{n}$ is the state vector, and $\mathbf{f}, \mathbf{G}$ are the non-linear functions. $\Gamma(t)$ is the Gaussian white-noise process of strength $\mathbf{Q}(t)$ and is used to account for the uncertainty in the dynamical model. Here, the initial state is described by the pdf p($\x_0$, $t_0$) and the pdf p($\x$, $t_k$) needs to be computed. For a general nonlinear dynamical system with Additive White Gaussian Noise (AWGN), the FPKE is given as:
% \begin{equation}\label{Eq:FPKEeqn}
% \begin{aligned}
% \frac{\partial p(\mathbf{x}, t)}{\partial t}=&-\frac{\partial p(\mathbf{x}, t)^T}{\partial \mathbf{x}} \mathbf{f}(\mathbf{x}, t)-p(\mathbf{x}, t) \cdot T r\left[\frac{\partial \mathbf{f}(\mathbf{x}, t)}{\partial \mathbf{x}}\right] +\frac{1}{2} \operatorname{Tr}\left[\mathbf{G} \mathbf{Q}_p(t) \mathbf{G}^T \frac{\partial^2 p(\mathbf{x}, t)}{\partial \mathbf{x} \partial \mathbf{x}^T}\right]
% \end{aligned}
% \end{equation}
\begin{equation}\label{Eq:FPKEeqn}
\begin{aligned}
\frac{\partial p(\mathbf{x}, t)}{\partial t} =\sum_{i=1}^n \frac{\partial\left(p f_i\right)}{\partial x_i}+\frac{1}{2} \sum_{i, j=1}^n \frac{\partial^2\left[\left(\mathbf{G Q G}^T\right)_{i j} p\right]}{\partial x_i \partial x_j}
\end{aligned}
\end{equation}
Notice that on the right-hand side, the first term represents the drift and the second term represents the diffusion of system response. The objective is to find the state pdf at time $t_k$, given by $p(\mathbf{x}, t_k)$ by solving the FPKE. Although the FPKE is linear in state pdf, it is a formidable equation to solve due to following reasons: 
\begin{enumerate}
    \item Positivity of the pdf: $p(\x, t) \geq 0 \quad \forall \mathbf{x}, t$
    \item Infinite Boundary Conditions of the pdf: $p(\pm \infty, t)=0$
    \item Normality of the pdf: $\int p(\mathbf{x}, t) d \mathbf{x}=1$
    \item Selection of appropriate Basis Functions
\end{enumerate}

Let us handle these constraints individually as we attempt to solve the FPKE. To address the non-negativity constraint, the following exponential form is assumed for the state pdf:
\begin{equation}\label{Eq:expBeta}
p(\mathbf{x}, t)=e^{\beta(\mathbf{x}, t)}
\end{equation}
where $\beta(\mathbf{x}, t)$ denotes the log-pdf. In an attempt to impose the infinite boundary conditions, the true pdf will be regularized by a weighting function $W(\mathbf{x}, t, \theta)$ that must satisfy the following properties: i)  $W(\mathbf{x}, t, \boldsymbol{\theta}) \geq 0 \quad \forall \mathbf{x}, t, \boldsymbol{\theta}$. ii) $W(-\infty, t, \theta)=W(t, \infty, \theta)=0 \quad \forall t, \theta$. Here $\boldsymbol{\theta}$ is a parameter vector based on the weighting function chosen. These two properties lead to the use of a probability density function as the regularizing function. Note that the weight function has a similar form to the true pdf given in \eqref{Eq:expBeta}. 
\begin{equation}\label{Eq:WeightFn}
W(\mathbf{x}, t, \boldsymbol{\theta})=e^{\beta_W(\mathbf{x}, t, \boldsymbol{\theta})}
\end{equation}
The weighting function is used to regularize the true pdf as follows:
\begin{equation}\label{Eq:p_approx}
p_A(\mathbf{x}, t)=p(\mathbf{x}, t) W(\mathbf{x}, t, \boldsymbol{\theta})=e^{\beta(\mathbf{x}, t)} e^{\beta_W(\mathbf{x}, t, \boldsymbol{\theta})}=e^{\left(\beta(\mathbf{x}, t)+\beta_W(\mathbf{x}, t, \boldsymbol{\theta})\right)}=e^{\beta_A(\mathbf{x}, t, \theta)}
\end{equation}
The log-pdf can now be calculated using \eqref{Eq:p_approx} with the positive and infinite boundary condition restrictions met. To handle the normality constraint, the pdf can be integrated over the domain of interest to find the normality constant. If the log-pdf is smooth with continuity and the existence of at least the first two derivatives, then $\beta_A(\mathbf{x}, t)$ can be represented as follows:
\begin{equation}\label{Eq:Beta_A_coeff_w_t}
\beta_A(\mathbf{x}, t)=\beta(\mathbf{x}, t)+\beta_W(\mathbf{x}, \boldsymbol{\theta})=\mathbf{c}^T(t) \Phi(\mathbf{x})+\mathbf{c}_W^T(\boldsymbol{\theta}) \Phi(\mathbf{x})
\end{equation}
where $\mathbf{c}(t) \in \mathbb{R}^m$ denotes a vector of time-varying coefficients of the $\log -\mathrm{pdf}, \mathbf{c}_W(\boldsymbol{\theta})$ is a vector of coefficients of the weighting function, and $\Phi(\mathbf{x}) \in \mathbb{R}^m$ is a vector of basis functions. 
The last constraint is the selection of the appropriate basis functions in finding the solution of FPKE. The fact that the number of spatial variables is equal to the state dimension is the main barrier to finding a numerical solution to the FPKE. As the state dimension increases, the number of coefficients corresponding to the dictionary of basis functions will increase rapidly. Therefore, selecting an appropriate basis function is crucial in computing an accurate solution. There are infinite choices to select the dictionary of basis function in approximating the log-pdf. These can be polynomial functions, B-spline functions, wavelet functions, etc. The polynomial functions, particularly monomials, are generally chosen as basis functions due to the agreement with the Principle of Maximum Entropy (PME) solution form and has been shown in our previous work \cite{jain2020computationally}. In this work, however, new basis functions exploiting the Hamiltonian structure are introduced along with monomials to capture and propagate the uncertainty accurately and efficiently.

To motivate for the Hamiltonian-type formulation, let us consider a stationary solution of the FPKE defined by an $n$ - dimensional dynamical system $\dot{\mathbf{x}}=\mathbf{f}(\mathbf{x}), n=2 N$. The first $N$ components correspond to position $(\mathbf{x})$ and the last $N$ to velocity $(\dot{\mathbf{x}})$.
Let $\mathbf{y}=\dot{\mathbf{x}}$ i.e.
\begin{equation} 
\begin{aligned}
\frac{d x_i}{d t}=f_i ; & \qquad  \frac{d y_i}{d t}=f_{i+N} ; & \quad   i=1, \cdots, N
\end{aligned}
\end{equation}
Lets define new generalized coordinates to define Hamiltonian as a function of configuration variable, $\mathbf{q}$ and corresponding conjugate momenta, $\mathbf{p}$, where $\mathbf{q} = \mathbf{x}$ and $\mathbf{p} = \mathbf{y}$. Hamiltonian is considered the system’s total energy (sum of kinetic energy and potential energy) and is conservative during the motion of the system :
% The Hamiltonian is given as : 
\begin{equation}\label{Hamiltonian_generalsys}
H(\mathbf{q}, \mathbf{p})=\sum_{i=1}^n \frac{p_i^2}{2 m_i}+V(\mathbf{q})
\end{equation}
Utilizing the Hamiltonian, the following relations can be obtained:
\begin{equation}
\begin{aligned}
&\frac{d \mathbf{x}}{d t}=\frac{d \mathbf{q}}{d t}=\frac{\partial H}{\partial \mathbf{p}} ; \qquad \frac{d \mathbf{y}}{d t}=\frac{d \mathbf{p}}{d t}=-\frac{\partial H}{\partial \mathbf{q}}
\end{aligned}
\end{equation}
% The FPKE given by \eqref{Eq:FPKEeqn} can be written for this system as:
The FPKE given in \eqref{Eq:FPKEeqn} can be solved by substituting $f_i$ in terms of Hamiltonian,
\begin{equation}
\begin{aligned}
\frac{\partial p(\mathbf{x}, t)}{\partial t} &=-\sum_{i=1}^n \frac{\partial\left(p f_i\right)}{\partial x_i}+\frac{1}{2} \sum_{i, j=1}^n \frac{\partial^2\left[\left(\mathbf{g Q g}^T\right)_{i j} p\right]}{\partial x_i \partial x_j} \\
&=-\sum_{i=1}^N \frac{\partial}{\partial x_i}\left(p \frac{d x_i}{d t}\right)+\frac{\partial}{\partial y_i}\left(p \frac{d y_i}{d t}\right)+\frac{1}{2} \sum_{i, j=1}^n \frac{\partial^2\left[\left(\mathbf{g Q g}^T\right)_{i j} p\right]}{\partial x_i \partial x_j} \\
&=-\sum_{i=1}^N\left(\frac{\partial}{\partial q_i}\left(p \frac{\partial H}{\partial p_i}\right)-\frac{\partial}{\partial p_i}\left(p \frac{\partial H}{\partial q_i}\right)\right)+\frac{1}{2} \sum_{i, j=1}^n \frac{\partial^2\left[(\mathbf{g Q g})_{i j} p\right]}{\partial x_i \partial x_j} \\
% &=-\sum_{i=1}^N\left(\frac{\partial p}{\partial q_i} \frac{\partial H}{\partial p_i}-\frac{\partial p}{\partial p_i} \frac{\partial H}{\partial q_i}+p \frac{\partial^2 H}{ \partial q_i \partial p_i}-p \frac{\partial^2 H}{ \partial p_i \partial q_i }\right) +\frac{1}{2} \sum_{i, j=1}^n \frac{\partial^2\left[(\mathbf{g Q g})_{i j} p\right]}{\partial x_i \partial x_j} \\
% \frac{\partial p(\mathbf{x}, t)}{\partial t} & = -\sum_{i=1}^N\left(\frac{\partial p}{\partial q_i} \frac{\partial H}{\partial p_i}-\frac{\partial p}{\partial p_i} \frac{\partial H}{\partial q_i} \right) +\frac{1}{2} \sum_{i, j=1}^n \frac{\partial^2\left[(\mathbf{g Q g})_{i j} p\right]}{\partial x_i \partial x_j} 
\end{aligned}
\end{equation}
\begin{equation}
\begin{aligned}
\frac{\partial p(\mathbf{x}, t)}{\partial t}  &= -\sum_{i=1}^N\left(\frac{\partial p}{\partial q_i} \frac{\partial H}{\partial p_i}-\frac{\partial p}{\partial p_i} \frac{\partial H}{\partial q_i}+p \frac{\partial^2 H}{ \partial q_i \partial p_i}-p \frac{\partial^2 H}{ \partial p_i \partial q_i }\right) +\frac{1}{2} \sum_{i, j=1}^n \frac{\partial^2\left[(\mathbf{g Q g})_{i j} p\right]}{\partial x_i \partial x_j} \\
\frac{\partial p(\mathbf{x}, t)}{\partial t} & = -\sum_{i=1}^N\left(\frac{\partial p}{\partial q_i} \frac{\partial H}{\partial p_i}-\frac{\partial p}{\partial p_i} \frac{\partial H}{\partial q_i} \right) +\frac{1}{2} \sum_{i, j=1}^n \frac{\partial^2\left[(\mathbf{g Q g})_{i j} p\right]}{\partial x_i \partial x_j} 
\end{aligned}
\end{equation}
Let $p(\mathbf{x}, t)=p\left(H\left(\x, \mathbf{y}\right)\right)=p\left(H\left(\mathbf{q}, \mathbf{p}\right)\right)$. Using the fact that for a stationary pdf $\frac{\partial p(\mathbf{x}, t)}{\partial t} =0$, the above equation can be written as:
\begin{equation}
0=-\sum_{i=1}^N\left(\frac{\partial p}{\partial H}\frac{\partial H}{\partial q_i} \frac{\partial H}{\partial p_i}-\frac{\partial p}{\partial H} \frac{\partial H}{\partial p_i}\frac{\partial H}{\partial q_i}\right) +\frac{1}{2} \sum_{i, j=1}^n \frac{\partial^2\left[(\mathbf{g Q g})_{i j} p\right]}{\partial x_i \partial x_j}  
\end{equation}
If $\mathbf{Q}=0$, any function of $H\left(\mathbf{q},\mathbf{p}\right)$ will serve for $p\left(H\left(\mathbf{q}, \mathbf{p}\right)\right)$ provided that it satisfies normality and boundary condition.
Otherwise, $p\left(H\left(\mathbf{q},\mathbf{p}\right)\right)$ is found by solving:
\begin{equation}
\frac{1}{2} \sum_{i, j=1}^n \frac{\partial^2\left[(\mathbf{g Q g})_{i j} p\right]}{\partial x_i \partial x_j} = 0 
\end{equation}

%to include Hamiltonian in the basis function dictionary of stationary pdf
To further get a feel about the stationary pdf dependence on Hamiltonian, let us take a numerical example of Duffing oscillator \cite{muscolino1997stationary}. The system dynamics are given as follows:
\begin{equation}\label{DO_EOMs}
\ddot{x}+\eta \dot{x}+\alpha x+\beta x^3= Q
\end{equation}
It can be written in state-space form as:  
\begin{equation}
\dot{x}_1=x_2 ; \qquad \dot{x}_2= Q -\eta x_2-\alpha x_1-\beta x_1^3
\end{equation}
The equations of motion can be written by defining Hamiltonian, $H=\frac{x_2^2}{2}+\frac{\alpha x_1^2}{2}+\frac{\beta x_1^4}{4}
$
\begin{equation}
\frac{d x_1}{d t}=\frac{\partial H}{\partial x_2} ; \qquad  \frac{d x_2}{d t}=-\frac{\partial H}{\partial x_1}-\eta \frac{\partial H}{\partial x_2}+Q
\end{equation}
The FPKE given by \eqref{Eq:FPKEeqn} can be written for this system as:
\begin{equation}
-\frac{\partial}{\partial x_1}\left(p \frac{\partial H}{\partial x_2}\right)+\frac{\partial}{\partial x_2}\left(p \frac{\partial H}{\partial x_1}\right)+\eta \frac{\partial}{\partial x_2}\left(p \frac{\partial H}{\partial x_2}\right)+\frac{1}{2} Q \frac{\partial^2 p}{\partial x_2^2}=0
\end{equation}
The first two terms cancel each other out, as was previously calculated. Additional calculations using the remaining terms result in:
\begin{equation}
\begin{aligned}
\eta \frac{\partial}{\partial x_2}\left(p \frac{\partial H}{\partial x_2}\right)+\frac{1}{2} Q \frac{\partial^2 p}{\partial x_2^2} &=0 \\
\eta p \frac{\partial H}{\partial x_2}+\frac{1}{2} Q \frac{\partial p}{\partial x_2} &=C \text { (constant) }
\end{aligned}
\end{equation}
Using $p(\pm \infty, t)=0 $
$$\eta p \frac{\partial H}{\partial x_2}+\frac{1}{2} Q \frac{\partial p}{\partial x_2}=0
$$
Further simplifications lead to:
\begin{equation}
\begin{aligned}
\eta p \frac{\partial H}{\partial x_2}+\frac{1}{2} Q \frac{\partial p}{\partial H} \frac{\partial H}{\partial x_2} &=0 \\
\eta p+\frac{1}{2} Q \frac{\partial p}{\partial H} &=0 \\
\end{aligned}
\end{equation}
\begin{equation}\label{Eq:DO_statPDF}
p =p_0 \exp \left(-\frac{2 \eta H}{Q} \right)  \\
\end{equation}

As shown above, the stationary pdf of a Duffing oscillator is a function of the Hamiltonian, and log-pdf of Duffing oscillator is directly proportionate to the Hamiltonian. Thus, this motivates us to include the Hamiltonians in the over-complete dictionary of basis functions along with monomials. Monomial basis functions are necessary to approximate the transient behavior of the log-pdf, while the Hamiltonian basis functions are necessary to approximate the stationary log-pdf. The complete basis function dictionary can be written as:
\begin{equation}
\Phi(\mathbf{x}_j) = [\phi(\mathbf{x}_j) , \phi(H)]
\end{equation}
where $\Phi(\mathbf{x}_j) \in \re^{m = m_x +m_h}$ is the dictionary of basis functions containing monomials and Hamiltonians, $\phi(\mathbf{x}_j) \in \re^{m_x}$ is the dictionary of basis functions containing monomials and $\phi(H) \in \re^{m_h}$ is the dictionary of basis functions containing different order of Hamiltonians. $\phi(H)$ can be written as:
\begin{equation}
\phi(H) = [H^1, H^2, H^3, \dots, H^{m_h}]
\end{equation}

After satisfying all necessary conditions and selecting appropriate basis functions, the next step is the computation of the coefficients from the linear constraint FPKE equation (derived in the next section). Various numerical techniques exist to solve the coefficients from this equality constraint equation. Techniques such as Galerkin-based methods \cite{beard1997galerkin} and least-squares can be applied to minimize the global residual error when a series approximation is employed for the log-pdf. However, these methods involve computing high-dimensional integrals and, thus, suffer from the curse of dimensionality. On the other hand, utilizing collocation-based methods to minimize the global residual error reduces the problem to the function evaluation at the given collocation/grid points. A recently developed Quadrature scheme known as Conjugate Unscented Transform (CUT) provides a minimal set of cubature points. Thus, CUT points are used as the collocation points to represent the domain accurately. Furthermore, a minimal expansion for the log-pdf approximation is desired. If all of the coefficients participate in capturing the pdf, it can lead to the overfitting of the data and consequently results in a high testing error. Thus, the sparse approximation technique will be utilized to judicially select the dominant coefficients from the over-complete dictionary of basis functions and will thus help in avoiding the overfitting of the data. 

% In the sparse-collocation method, the first step is to approximate the log-pdf using a polynomial-based expansion. Then a set of collocation points needs to be generated to represent the pdf domain perfectly. A recently developed Quadrature scheme known as Conjugate Unscented Transform (CUT) provides a minimal set of cubature points. Thus, CUT points are used as the collocation points to represent the domain accurately. Further, a minimal expansion for the log-pdf approximation is desired. However, for polynomial basis functions, the number of required polynomials is combinatorial and thus proliferates for the desired order expansion. For th order polynomials in ��-dimensional space, the required number of basis functions is ... . This growth quickly outpaces the CUT-generated collocation points, notably for higher dimensional systems. In order to obtain a minimal expression for the value function approximation without affecting the number of collocation points, this research utilizes a sequential ��1-norm minimization routine to optimally select the basis functions from the over-complete set�

In the following section, a sparse-collocation technique is described to solve the FPKE in a computationally attractive manner. The proposed approach exploits the recent advances in non-product cubature methods in conjunction with carefully selecting the basis functions and sparse approximation tools to alleviate the curse of dimensionality to some extent.

\subsection{Development of Collocation Equations} \label{CollocationEquations}
In this section, the numerical methodology for solving the FPKE using the collocation method is briefly mentioned. The full methodology has been derived in \cite{mercurio2016conjugate}. As mentioned in the last section, the non-negativity and infinite boundary conditions are handled using the log-pdf and a weight function, as seen in \eqref{Eq:p_approx}. To develop the collocation equations using FPKE, \eqref{Eq:p_approx} can be substituted into \eqref{Eq:FPKEeqn} and the weighting function is assumed to be a constant Gaussian kernel. Further, all-time derivatives can be expanded as finite differences and the method of weighted residuals can be used to derive the set of equations governing each of the $m$ coefficients:

\begin{equation}
\begin{aligned}
e(\mathbf{x}, t) &=\Phi^T(\mathbf{x})\left(\frac{\mathbf{c}\left(t_{k+1}\right)-\mathbf{c}\left(t_k\right)}{\Delta t}\right)+\operatorname{Tr}\left[\frac{\partial \mathbf{f}\left(t_k, \mathbf{x}\right)}{\partial \mathbf{x}}\right] +\mathbf{f}^T\left(t_k, \mathbf{x}\right)\left[\frac{\partial \Phi(\mathbf{x})^T}{\partial \mathbf{x}}\left(\mathbf{c}\left(t_k\right)+\mathbf{c}_W(\boldsymbol{\theta})\right)\right] \\
&-\frac{1}{2} \operatorname{Tr}\left[\mathbf { g } ( t _ { k } ) \mathbf { Q } ( t _ { k } ) \mathbf { g } ( t _ { k } ) ^ { T } \left(\left(\mathbf{c}\left(t_k\right)+\mathbf{c}_W(\boldsymbol{\theta})\right)^T \frac{\partial \Phi(\mathbf{x})}{\partial \mathbf{x}} \frac{\partial \Phi(\mathbf{x})}{\partial \mathbf{x}^T}\left(\mathbf{c}\left(t_k\right)+\mathbf{c}_W(\boldsymbol{\theta})\right)\right.\right.\\
&\left.\left.+\frac{\partial^2}{\partial \mathbf{x} \partial \mathbf{x}^T}\left[\left(\mathbf{c}\left(t_k\right)+\mathbf{c}_W(\boldsymbol{\theta})\right)^T \Phi(\mathbf{x})\right]\right)\right] -\frac{\partial \mathbf{c}_W(\boldsymbol{\theta})^T}{\partial \boldsymbol{\theta}}\left(\frac{\boldsymbol{\theta}\left(t_{k+1}\right)-\boldsymbol{\theta}\left(t_k\right)}{\Delta t}\right)
\end{aligned}
\end{equation}

In the collocation method, the residual error is projected onto a series of delta functions centered at chosen collocation points resulting in a residual error being zero at the collocation points. The selection of the collocation points is crucial in obtaining a well-conditioned system of equations for the unknown coefficients. Assuming there are total $N$ collocation points, leading to a system of $N$ equations in $m$ unknowns to exactly solve the FPKE at prescribed points, $\mathbf{x}_i$ :
\begin{equation}
\int e(\mathbf{x}, t) \delta\left(\mathbf{x}-\mathbf{x}_i\right) d \mathbf{x}=0 \rightarrow e\left(\mathbf{x}, t_i\right)=0, \quad i=1,2, \ldots N
\end{equation}
where $\mathbf{x}_i$ are the chosen collocation points. This leads to the following system of equations for the unknown coefficients:
% \begin{equation}\label{Eq:ck+1_ck}
% \mathbf{A} \frac{\mathbf{c}\left(t_{k+1}\right)-\mathbf{c}\left(t_k\right)}{\Delta t}+\mathbf{b}=0
% \end{equation}
% where the $j^{t h}$ row of the matrix $\mathbf{A}$ and vector $\mathbf{b}$ are given as: 
% \begin{equation}
% \mathbf{A}_j=\Phi^T\left(\mathbf{x}_j\right)
% \end{equation}
% \begin{equation}
% \begin{aligned}
% &\mathbf{b}_j=\left[\operatorname{Tr}\left(\frac{\partial \mathbf{f}\left(t_k, \mathbf{x}\right)}{\partial \mathbf{x}}\right)+\mathbf{f}^T\left(t_k, \mathbf{x}\right)\left(\frac{\partial \Phi(\mathbf{x})^T}{\partial \mathbf{x}}\left(\mathbf{c}\left(t_k\right)+\mathbf{c}_W\right)\right)\right] \\
% &-\frac{1}{2} \operatorname{Tr}\left[\mathbf { g } ( t _ { k } ) \mathbf { Q } ( t _ { k } ) \mathbf { g } ^ { T } ( t _ { k } ) \left(\left(\mathbf{c}\left(t_k\right)+\mathbf{c}_W\right)^T \frac{\partial \Phi(\mathbf{x})}{\partial \mathbf{x}} \frac{\partial \Phi(\mathbf{x})}{\partial \mathbf{x}^T}\left(\mathbf{c}\left(t_k\right)+\mathbf{c}_W\right)\right.\right.\\
% &\left.+\frac{\partial^2}{\partial \mathbf{x} \partial \mathbf{x}^T}\left[\left(\mathbf{c}\left(t_k\right)+\mathbf{c}_W\right)^T \Phi(\mathbf{x})\right]\right)  \left.-\frac{\partial \mathbf{c}_W(\boldsymbol{\theta})^T}{\partial \boldsymbol{\theta}}\left(\frac{\boldsymbol{\theta}\left(t_{k+1}\right)-\boldsymbol{\theta}\left(t_k\right)}{\Delta t}\right)\right]_{\mathbf{x}=\mathbf{x}_j}
% \end{aligned}
% \end{equation}
\begin{equation}\label{Eq:ck+1_ck}
\mathbf{A} \mathbf{c}_{k+1} +\mathbf{b}=0
\end{equation}
% where the $j^{t h}$ row of the matrix $\mathbf{A}$ and vector $\mathbf{b}$ are given as: 
% \begin{equation}
% \mathbf{A}_j=\Phi^T\left(\mathbf{x}_j\right)
% \end{equation}
% \begin{equation}
% \begin{aligned}
% & \mathbf{A}_j=\Phi^T\left(\mathbf{x}_j\right) \\
% & \mathbf{b}_j=  - \Phi^T\left(\mathbf{x}_j\right) \mathbf{c}_k   + \left[\operatorname{Tr}\left(\frac{\partial \mathbf{f}\left(t_k, \mathbf{x}\right)}{\partial \mathbf{x}}\right)+\mathbf{f}^T\left(t_k, \mathbf{x}\right)\left(\frac{\partial \Phi(\mathbf{x})^T}{\partial \mathbf{x}}\left(\mathbf{c}_k+\mathbf{c}_W\right)\right) \\
% &-\frac{1}{2} \operatorname{Tr}\left[\mathbf { g } ( t _ { k } ) \mathbf { Q } ( t _ { k } ) \mathbf { g } ^ { T } ( t _ { k } ) \left(\left(\mathbf{c}_k+\mathbf{c}_W\right)^T \frac{\partial \Phi(\mathbf{x})}{\partial \mathbf{x}} \frac{\partial \Phi(\mathbf{x})}{\partial \mathbf{x}^T}\left(\mathbf{c}_k+\mathbf{c}_W\right)\right.\right. \left.+\frac{\partial^2}{\partial \mathbf{x} \partial \mathbf{x}^T}\left[\left(\mathbf{c}_k+\mathbf{c}_W\right)^T \Phi(\mathbf{x})\right]\right)  \left.\right]_{\mathbf{x}=\mathbf{x}_j} \right] \Delta t
% \end{aligned}
% \end{equation}
\begin{equation}
\begin{gathered}
\mathbf{A}^j_i=\boldsymbol \Phi_j^T\left(\mathbf{x}_i\right) \\
\\
\mathbf{b}_i=  \boldsymbol \Phi^T\left(\mathbf{x}_i\right) \mathbf{c}_k   - \Biggl[\operatorname{Tr}\left(\frac{\partial \mathbf{f}\left(t_k, \mathbf{x}\right)}{\partial \mathbf{x}}\right)+\mathbf{f}^T\left(t_k, \mathbf{x}\right)\left(\frac{\partial \boldsymbol \Phi^T(\mathbf{x})}{\partial \mathbf{x}}\left(\mathbf{c}_k+\mathbf{c}_W\right)\right) \\
-\frac{1}{2} \operatorname{Tr} \Biggl\{ \mathbf{g} (t_{k}) \mathbf{Q} (t_{k} ) \mathbf { g } ^ { T } ( t _ { k } ) \Biggl( ( \mathbf{c}_k+\mathbf{c}_W )^T \frac{\partial \boldsymbol \Phi(\mathbf{x})}{\partial \mathbf{x}} \frac{\partial \boldsymbol \Phi(\mathbf{x})}{\partial \mathbf{x}^T} (\mathbf{c}_k+\mathbf{c}_W )  \\  
+\frac{\partial^2}{\partial \mathbf{x} \partial \mathbf{x}^T}\left[ (\mathbf{c}_k+\mathbf{c}_W )^T \boldsymbol \Phi(\mathbf{x})\right] \Biggr)   \Biggr\} \Biggr]_{\mathbf{x}=\mathbf{x}_i} \Delta t,   \qquad i = 1 , 2, \dots ,N \qquad j = 1,2, \dots ,m
\end{gathered}
\end{equation}

Notice that the numerical solution of \eqref{Eq:ck+1_ck} can be computed using the least-squares method, which determines the best fit solution for the given $N$ collocation points resulting in $m $ coefficients. The least-squares method utilizes most of these $m $ coefficients from the dictionary of basis functions and thus provides the smallest possible two-norm error. As a result, this method tends to overfit the training data (collocation points) and can yield a high norm error on the testing data. Therefore, an alternative method utilizing $l_1$-norm approximation is implemented to provide the minimum possible number of coefficients required to represent the entire domain accurately. In addition, the objective is to find the minimum possible number of collocation points that can accurately represent the entire domain and prevent the infeasibility of the solution in higher dimensions due to the curse of dimensionality. %The $m$ coefficients computed with this method provide the minimum norm error by utilizing most of the basis function coefficients.

\subsection{Selection of Collocation Points}
As mentioned in the previous section, the collocation points should be selected so that a well-conditioned system of equations is obtained. The traditional Monte Carlo method is easier to implement but has a slow convergence rate. The alternative method employs numerical sampling techniques, namely quadrature methods. Many methods exist for generating quadrature points, with the simplest method being Gaussian quadrature. The Gaussian quadrature methods provide a minimal number of quadrature points to integrate the polynomial function in a one-dimensional (1-D) space. However, one needs to take a tensor product of these 1-D quadrature points to evaluate multidimensional integrals in $d$-D space. This results in a total of $q^d$ quadrature points, where $q$ is the number of quadrature points in 1-D. Thus, the number of these points increases exponentially as the system's dimension increases\cite{numintsparsegrid,strGQF}. 

An efficient numerical sampling method known as Conjugate Unscented Transformation (CUT) method is chosen as collocation points to prevent this exponential growth. This method exploits the structure of the domain to choose specially defined axes on which quadrature points are selected. These CUT points are developed using minimal cubature rules and offer similar orders of accuracy as Gaussian quadratures with fewer required points\cite{venkat,adurthi2018conjugate}. More details about the CUT methodology and its comparison with conventional quadrature rules can be found in Ref.~\cite{venkat_thesis,cutgaussgnc,adurthi2015conjugate,venkat_jgcd}.

Therefore, these CUT collocation points provide the locations at which the polynomial expansion satisfies the FPKE. Thus, it is crucial to choose the parameters of the CUT method so that the sigma points accurately sample the true pdf. Further, to approximate the polynomial coefficients of the log-pdf, the collocation points or the basis function must be transformed to a specific range. For example, the collocation points can be mapped uniformly into $\pm 1$ or the domain of zero mean and identity covariance. This is needed due to the fact that optimization algorithms are sensitive to magnitude and range variations across each dimension of the system, which, if not handled properly, can lead to numerical issues. Assuming the global domain is known a-priori, a linear transformation can be used to map the global domain of interest into a local domain within a hypercube as:
\begin{equation} \label{Eq:y=T(x+B)}
\mathbf{y}=\mathbf{T_0}\left(\mathbf{x}+\mathbf{B}_0\right)
\end{equation}
% Instead of mapping the polynomial basis to the local domain, it is simpler to rewrite the FPKE in the local space and map only the system dynamics and Jacobian matrix via \eqref{Eq:y=T(x+B)} as:
% \begin{equation}
% \begin{gathered}
% \dot{\mathbf{y}}=\mathbf{T}_0 \dot{\mathbf{x}}=\mathbf{T}_0\left[\mathbf{f}\left(\mathbf{x}, t=\mathbf{T}_0^{-1} \mathbf{y}-\mathbf{B}_0\right)+\mathbf{G}\left(\mathbf{x}, t=\mathbf{T}_0^{-1} \mathbf{y}-\mathbf{B}_0\right) \Gamma(t)\right]=\overline{\mathbf{f}}(\mathbf{y}, t)+\overline{\mathbf{G}}(\mathbf{y}, t) \Gamma(t) \\
% \frac{d \overline{\mathbf{f}}(t, \mathbf{y})}{d \mathbf{y}}=\left[\mathbf{T}_0 \frac{d \mathbf{f}(\mathbf{x}, t)}{d \mathbf{x}} \mathbf{T}_0^{-1}\right]_{\mathbf{x}=\mathbf{T}_0^{-1} \mathbf{y}-\mathbf{B}_0}
% \end{gathered}
% \end{equation}
The transformation above results in the following revised set of ODEs:
\begin{equation}\label{Eq:ck+1_ck}
\mathbf{A} \mathbf{c}_{k+1} +\mathbf{b}=0
\end{equation}
where the $j^{t h}$ row of the matrix $\mathbf{A}$ and vector $\mathbf{b}$ are given as: 
\begin{equation}
\mathbf{A}_j=\Phi^T\left(\mathbf{y}_j\right)
\end{equation}
% \begin{equation}\label{Eq:bj(y)}
% \begin{aligned}
% &\mathbf{b}_j=  - \Phi^T\left(\mathbf{y}_j\right) \mathbf{c}_k   + \left[\operatorname{Tr}\left(\frac{\partial \mathbf{f}\left(t_k, \mathbf{y}\right)}{\partial \mathbf{y}}\right)+\mathbf{f}^T\left(t_k, \mathbf{y}\right)\left(\frac{\partial \Phi(\mathbf{y})^T}{\partial \mathbf{y}}\left(\mathbf{c}_k+\mathbf{c}_W\right)\right) \\
% &-\frac{1}{2} \operatorname{Tr}\left[\mathbf { g } ( t _ { k } ) \mathbf { Q } ( t _ { k } ) \mathbf { g } ^ { T } ( t _ { k } ) \left(\left(\mathbf{c}_k+\mathbf{c}_W\right)^T \frac{\partial \Phi(\mathbf{y})}{\partial \mathbf{y}} \frac{\partial \Phi(\mathbf{y})}{\partial \mathbf{y}^T}\left(\mathbf{c}_k+\mathbf{c}_W\right)\right.\right. \left.+\frac{\partial^2}{\partial \mathbf{y} \partial \mathbf{y}^T}\left[\left(\mathbf{c}_k+\mathbf{c}_W\right)^T \Phi(\mathbf{y})\right]\right)  \left. \right]_{\mathbf{y}=\mathbf{y}_j} \right] \Delta t
% \end{aligned}
% \end{equation}
\begin{equation}\label{Eq:bj(y)}
\begin{gathered}
\mathbf{b}_i =  - \Phi^T\left(\mathbf{y}_i\right) \mathbf{c}_k   + \Biggl[\operatorname{Tr}\left(\frac{\partial \mathbf{f}\left(t_k, \mathbf{y}\right)}{\partial \mathbf{y}}\right)+\mathbf{f}^T\left(t_k, \mathbf{y}\right)\left(\frac{\partial \Phi^T(\mathbf{y})}{\partial \mathbf{y}}\left(\mathbf{c}_k+\mathbf{c}_W\right)\right) \\
-\frac{1}{2} \operatorname{Tr} \Biggl\{ \mathbf { g } ( t _ { k } ) \mathbf { Q } ( t _ { k } ) \mathbf { g } ^ { T } ( t _ { k } ) \Biggl( ( \mathbf{c}_k+\mathbf{c}_W )^T \frac{\partial \Phi(\mathbf{y})}{\partial \mathbf{y}} \frac{\partial \Phi(\mathbf{y})}{\partial \mathbf{y}^T} (\mathbf{c}_k+\mathbf{c}_W ) \\  
+\frac{\partial^2}{\partial \mathbf{y} \partial \mathbf{y}^T}\left[ (\mathbf{c}_k+\mathbf{c}_W )^T \Phi(\mathbf{y})\right] \Biggr)   \Biggr\} \Biggr]_{\mathbf{y}=\mathbf{y}_i} \Delta t ,   \qquad i = 1 , 2, \dots ,N \qquad j = 1,2, \dots ,m
\end{gathered}
\end{equation}
% \begin{equation}
% \begin{gathered}
% \mathbf{A}^j_i=\boldsymbol \Phi_j^T\left(\mathbf{x}_i\right) \\
% \\
% \mathbf{b}_i=  \boldsymbol \Phi^T\left(\mathbf{x}_i\right) \mathbf{c}_k   - \Biggl[\operatorname{Tr}\left(\frac{\partial \mathbf{f}\left(t_k, \mathbf{x}\right)}{\partial \mathbf{x}}\right)+\mathbf{f}^T\left(t_k, \mathbf{x}\right)\left(\frac{\partial \boldsymbol \Phi^T(\mathbf{x})}{\partial \mathbf{x}}\left(\mathbf{c}_k+\mathbf{c}_W\right)\right) \\
% -\frac{1}{2} \operatorname{Tr} \Biggl\{ \mathbf{g} (t_{k}) \mathbf{Q} (t_{k} ) \mathbf { g } ^ { T } ( t _ { k } ) \Biggl( ( \mathbf{c}_k+\mathbf{c}_W )^T \frac{\partial \boldsymbol \Phi(\mathbf{x})}{\partial \mathbf{x}} \frac{\partial \boldsymbol \Phi(\mathbf{x})}{\partial \mathbf{x}^T} (\mathbf{c}_k+\mathbf{c}_W )  \\  
% +\frac{\partial^2}{\partial \mathbf{x} \partial \mathbf{x}^T}\left[ (\mathbf{c}_k+\mathbf{c}_W )^T \boldsymbol \Phi(\mathbf{x})\right] \Biggr)   \Biggr\} \Biggr]_{\mathbf{x}=\mathbf{x}_i} \Delta t,   \qquad i = 1 , 2, \dots ,N \qquad j = 1,2, \dots ,m
% \end{gathered}
% \end{equation}
The system of ODEs in \eqref{Eq:ck+1_ck} yields the coefficients of the log-pdf expansion in the local space, which can be mapped back to the global space under the following transformation of variables \cite{jazwinski2007stochastic}:
\begin{equation}
p(\mathbf{x}, t)=p \left(\mathbf{y}=\mathbf{T}_0 (\mathbf{x}+\mathbf{B}_0), t \right)\left|\frac{\partial \mathbf{y}}{\partial \mathbf{x}}\right|
\end{equation}
Furthermore, the initial coefficients are chosen in such a way that the initial pdf approximation corresponds to the true initial condition pdf. With the collocation points and equations obtained, the next objective of this research, the optimal selection of the basis functions, can be investigated. In the log-pdf approximation, an increase in the number of collocation points beyond the number of basis function terms necessitates an increase in the number of basis functions. Due to the combinatorial expansion of polynomial basis functions, the number of required basis functions quickly exceeds the collocation points. A further increase in the number of collocation points would exceed the number of basis functions, necessitating another increase in the number of basis functions. This process would repeat itself, and the computational load would quickly become unmanageable, even for lower-dimensional systems.

To avoid the back-and-forth increase in collocation points and basis functions, a minimal expansion for the log-pdf is sought by exploiting the fact that there are fewer collocation points than expansion terms. As a result, a methodology is required to determine the basis functions included in the series expansion without modifying the chosen collocation points.

\subsection{Selection of Optimal Coefficients}

\begin{algorithm}
\caption{Collocation-Based Solution of the Fokker-Planck-Kolmogorov Equation \label{algo:FPKE_algo}}
  \begin{algorithmic}[1] 
    \INPUT $\mathbf{f}(\mathbf{x}), \mathbf{g}(\mathbf{x}), m$ basis $\phi(\mathbf{x})$, initial values of coefficients $\mathbf{c}\left(t_0\right)$, weight function coefficients $\mathbf{c}_W$, discretized time vector $t$, time step $\Delta t$, collocation points $X_i$ $i=1,2, \ldots, N$, initial weight matrix $\mathbf{K}_0$, and weight update parameter $\epsilon$.
    \OUTPUT $\beta(\mathbf{x}, t)$.
    \STATE $\mathbf{K}=\mathbf{K}_0$.
    \STATE  Compute matrix $\mathbf{A}$.
    \FOR{t=0,t $< t_f$, k=k+1}
        \STATE Compute $\mathbf{b}$ using $\mathbf{c}_k$ 
        \STATE $\mathbf{c}_{k+1}^*=$ WeightedOpt$(\mathbf{A}, \mathbf{b},\mathbf{W},\mathbf{c}_{{k+1}_{l_2}}, \Delta_{s},\alpha, \epsilon_{ls}, \eta  )$
        \STATE $\mathbf{c}_{k+1}=$ ComputeRScoefficients($\mathbf{c}_{k+1}^*, \mathbf{A}, \mathbf{b}, \delta_{rs} )$
        \STATE $\beta\left(\mathbf{x}, t_{k+1}\right) = \left(\mathbf{c}_{k+1}+\mathbf{c}_W\right)^T \Phi(\mathbf{x})$
    \ENDFOR
  \end{algorithmic}
\end{algorithm}

As described earlier, the solution of \eqref{Eq:ck+1_ck} can be calculated by minimization the two-norm error, which aims to find the best-fit solution for the given collocation points (training data). The optimal value of coefficients $\mathbf{c}_{{k+1}_{l_2}}$ can be found by solving the weighted $l_2$-norm minimization:
\begin{equation}
\begin{aligned}
\mathbf{c}_{{k+1}_{l_2}} = & \min_{\mathbf{c}_{{k+1}}}  \left\|\mathbf{W}( \mathbf{A} \mathbf{c}_{k+1} + \mathbf{b} ) \right\|_{2} \\
\end{aligned}
\end{equation}
where weight matrix $\mathbf{W}$ can be chosen appropriately. For this study, the diagonal matrix corresponding to the CUT point weights is utilized. The $l_2$-norm coefficients are computed using:
\begin{equation}\label{WeightedLS}
\begin{aligned}
\mathbf{c}_{{k+1}_{l_2}} =& \mathbf{A}^{\dagger}  \mathbf{b}
\end{aligned}
\end{equation}
where $\mathbf{A}^{\dagger} $ represents the pseudo-inverse, and its value can be computed depending on whether the system is under-determined (m $>$ N) or over-determined (m $<$ N). If the systems of equation is over-determined, $\mathbf{A}^{\dagger}  = \left(\mathbf{A}^{T} \mathbf{W} \mathbf{A}\right)^{-1} \mathbf{A}^{T} \mathbf{W} $ and if the system is under-determined, $\mathbf{A}^{\dagger}  = \mathbf{A}^{T} \mathbf{W} \left(\mathbf{A}^{T} \mathbf{W} \mathbf{A}\right)^{-1} $. %The minimum weighted two-norm error $\epsilon_{ls}$ can additionally be computed utilizing the LS coefficients $\mathbf{c}_{{k+1}_{l_2}} $, which can then serve as the criterion for the two-norm error in sparse solution $ \epsilon_{ls}  =  \left\|\mathbf{W}( \mathbf{A}   \mathbf{c}_{{k+1}_{l_2}}   - \mathbf{b} ) \right\|_{2} $. 

As previously stated, $\mathbf{c}_{{k+1}_{l_2}}$ is known to pick all possible coefficients from the dictionary of basis functions and is therefore not sparse. This research seeks a minimal polynomial expansion that guarantees sparsity for the log-pdf. Therefore, a weighted $l_1$-norm optimization problem is proposed to select the minimum possible coefficients from an over-extensive dictionary of basis functions. In lieu of the equality constraint of \eqref{Eq:ck+1_ck}, this optimization problem considered bounded two-norm error as a soft inequality constraint. This allows sparse coefficients $\mathbf{c}_{{k+1}_s}$ to trade sparsity for approximation error, providing a more flexible option. Thus, the optimization problem can be stated in terms of the coefficients $\mathbf{c}_{{k+1}_s}$:

\begin{equation}
\begin{aligned}
&\min _{\mathbf{c}_{k+1}}\left\|\mathbf{K} \mathbf{c}_{k+1}\right\|_{1}\\
\text { subject to: } \\
& \left\|\mathbf{W}( \mathbf{A} \mathbf{c}_{k+1} + \mathbf{b} ) \right\|_{2}  \leq \epsilon 
\end{aligned}
\end{equation}
where $\mathbf{K}$ is the coefficient weight matrix, and $\epsilon$ represents the soft constraint of the two-norm error. The initial value of $\mathbf{K}$ can be set to 1 or based on any \textit{a-priori} knowledge such as penalizing the higher-order coefficients or using any information gained from $l_2$-norm solution. For subsequent iterations, $\mathbf{K}$ can be updated to penalize the non-dominant coefficients using the solution of the previous iteration $ \mathbf{c}^-_{{k+1}}$: 
$$\mathbf{K}= \frac{1} {(\mathbf{c}^-_{{k+1}} + \eta)}$$ 
where $\eta$ is a small number in order to prevent division by zero. The value of $\epsilon$ is chosen to provide a more flexible solution to the optimization problem, thereby balancing approximation error with sparsity. The algorithm \ref{SparseSolution_algo} describes the optimal procedure to select the minimal polynomial expansion of log-pdf. The user-defined parameter $\Delta_{s}$ represents the desired difference between the two norms of consecutive coefficients. The algorithm eventually arrives at an optimal solution $\mathbf{c}_{{k+1}_s}$, when the difference between the two norms of consecutive coefficient values, $\delta$, is less than $\Delta_{s}$.

% from $2\%-10\%$ from the norm error from the least-squares solution. 
%The coefficient matrix $\mathbf{K}$ is initialized as ones and updated as shown in Algorithm \ref{SparseSolution_algo} to find the optimal sparse coefficient $\mathbf{c}_{k+1}^*$. The parameter values used in this algorithm are assumed as $\eta = 1e-4$ and $\Delta_{s} = 1e-4 $, where $\eta$ is a small number used in the denominator to update weight matrix $\mathbf{K}$ so that it doesn't go to infinity when a zero coefficient value is encountered and $\Delta_{s}$ is the minimum two-norm difference of the consecutive coefficient value required for this algorithm to exit. The algorithm will converge to optimal solution $\mathbf{c}_{k+1}^*$, when 2-norm difference of the consecutive coefficient values, $\delta$ is less than $\Delta_{s}$. 

\begin{algorithm}
\caption{Iterative weighted $l_1$-norm optimization: $\mathbf{c}_{{k+1}_s} $= WeightedOpt($\mathbf{A}, \mathbf{b},\mathbf{W},\mathbf{c}_{{k+1}_{l_2}}, \Delta_{s},\alpha, \epsilon , \eta  )$  \label{SparseSolution_algo}}
  \begin{algorithmic}[1] 
    \INPUT $\mathbf{A}, \mathbf{b},\mathbf{W},\mathbf{c}_{{k+1}_{l_2}}, \Delta_{s}, \alpha, \epsilon , \eta$
    \OUTPUT $\mathbf{c}_{{k+1}_s}$
    \STATE \textbf{Initialization}   $\mathbf{K}= \frac{1} {(\mathbf{c}_{{k+1}_{l_2}} + \eta)}, \delta=1  $ 
    \STATE  compute $\mathbf{c}_{k+1}^- =  \displaystyle \min _{\mathbf{c}_{k+1}}\left\|\mathbf{K} \mathbf{c}_{k+1}\right\|_{1} $  
    
     $\quad \qquad  \text{subject to:} \quad \left\|\mathbf{W}( \mathbf{A} \mathbf{c}_{k+1} + \mathbf{b} ) \right\|_{2}  \leq \epsilon$ 
    \WHILE{$\delta \ge \Delta_{s}$}
      \STATE Update $\mathbf{K}= \frac{1}{(\mathbf{c}_{k+1}^- + \eta)}$, to find $\mathbf{c}_{k+1}^+  =  \displaystyle \min _{\mathbf{c}_{k+1}}\left\|\mathbf{K} \mathbf{c}_{k+1}\right\|_{1} $  
      
      $\quad \quad\quad \quad\quad \qquad \quad \qquad\quad \qquad  \text{subject to:} \quad \left\|\mathbf{W}( \mathbf{A} \mathbf{c}_{k+1} + \mathbf{b} ) \right\|_{2}  \leq \epsilon$
      \STATE  Compute $\delta = \| \mathbf{c}_{k+1}^+ - \mathbf{c}_{k+1}^-\|_2 $
      \STATE $\mathbf{c}_{k+1}^- = \mathbf{c}_{k+1}^+ $ 
    \ENDWHILE
    \STATE $\mathbf{c}_{{k+1}_s} = \mathbf{c}_{k+1}^-$
  \end{algorithmic}
\end{algorithm}

After obtaining the sparse coefficients $\mathbf{c}_{{k+1}_s}$, it is possible to separate the dominant coefficients from the non-dominant coefficients by choosing a user-defined coefficient threshold $\delta_{rs}$. In addition, the non-dominant coefficients can be ignored for computational purposes by substituting zero. Therefore, a new minimal representation of the basis functions $ \mathbf{A}_{rs}\in \mathbb{R}^{N\times m _r} $ corresponding to the $m _r$ non-zero coefficients can be constructed. Thus, the reduced sparse (RS) coefficients, $\mathbf{c}_{k+1} \in \mathbb{R}^{m _r}$, can be calculated by minimizing the $l_2$-norm error using the dominant basis functions:
%This can provide a minimal representation of polynomial coefficients along with the minimal basis functions and the log-pdf. A new reduced sparse solution $\mathbf{c}_{k+1}$ can then be computed using the weighted least square solution with the minimal basis functions representation, $\mathbf{A}_{rs}$, minimal log-pdf, $\mathbf{b}_{red}$ and Weight matrix, $\mathbf{W}_{red}$
\begin{equation}
\begin{aligned}
\mathbf{c}_{{k+1}} =& \mathbf{A}_{rs}^{\dagger}  \mathbf{b}
\end{aligned}
\end{equation}
where $\mathbf{A}_{rs}^{\dagger}$ is the pseudo-inverse, and the solution can be found in a manner analogous to \eqref{WeightedLS}. This minimal representation $\mathbf{c}_{{k+1}}$ is employed to compute the state-pdf and the marginalized state pdf.

Algorithm \ref{algo:FPKE_algo} illustrates the complete algorithm for calculating a minimal polynomial expression for log-pdf till the final time $t_f$. The $l_1$-norm minimization issue addressed in Algorithm \ref{SparseSolution_algo} is convex and must be invoked at each time instant to find the time-varying state pdf. Various efficient convex optimization solvers, including Sedumi, SDPT3, Gurobi, MOSEK, and GLPK, are available in the CVX programming language to find the optimal solution\cite{boyd2004convex,cvx}. While finding a sparse solution, it must be kept in mind that a single solver may not be able to provide the optimal solution, as each solver has different capabilities and performance levels. Since numerical methods for convex optimization are not exact, the results are computed within a numerical precision or tolerance that has been predefined. The best precision level at which CVX considers a model to be fully solved is $O(1e-8)$, and it may not be able to compute the results for the lower numerical tolerance. The solver may throw an Inaccurate/Solved solution flag if this occurs. In this case, it is advised to replace the solver or reduce the required precision tolerance. %As this work requires a time-varying solution, the Algorithm \ref{SparseSolution_algo} must be invoked at each time instant.   %Further, an optimal trade-off curve will be plotted by solving a family of optimization problems and varying the constraint in the numerical simulations. 

%% file: src/03_numeric.tex
In this section, numerical examples are demonstrated to show the applicability of the proposed approach. A nonlinear oscillator is examined as a benchmark problem to evaluate the precision of the suggested methodology since an analytical solution to stationary FPKE is known. Additionally, an orbit maneuver problem corresponding to satellite motion in a Low-Earth-Orbit (LEO) is considered to demonstrate the extensibility to the systems of higher dimensionality.

\subsection{Duffing Oscillator}

\begin{figure}
\centering
\subfigure[t = 2 s]{\label{fig:DOcaseII_PDFa} \includegraphics[width=0.23\textwidth]{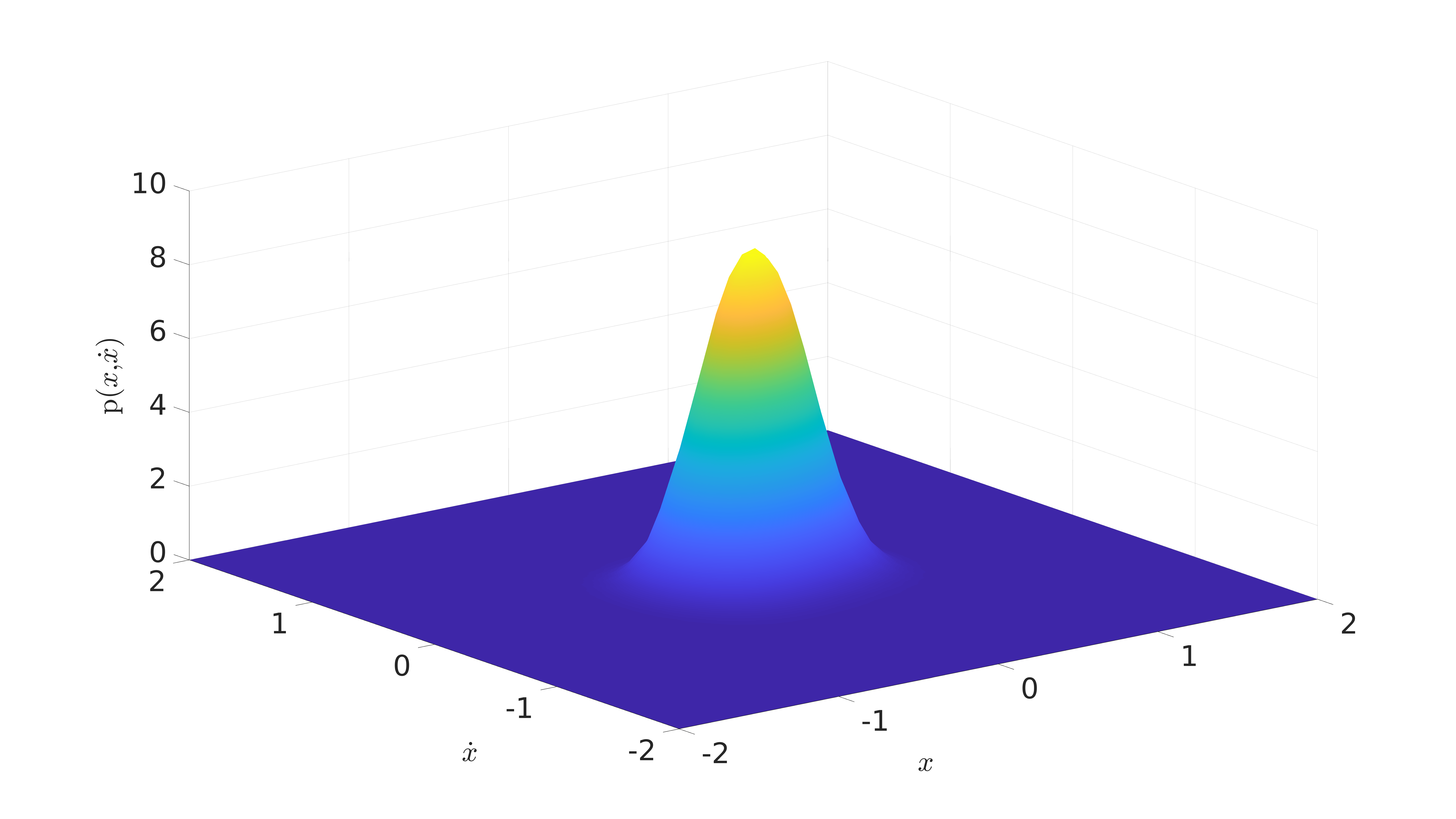}}
\subfigure[t = 5 s]{ \label{fig:DOcaseII_PDFb} \includegraphics[width=0.23\textwidth]{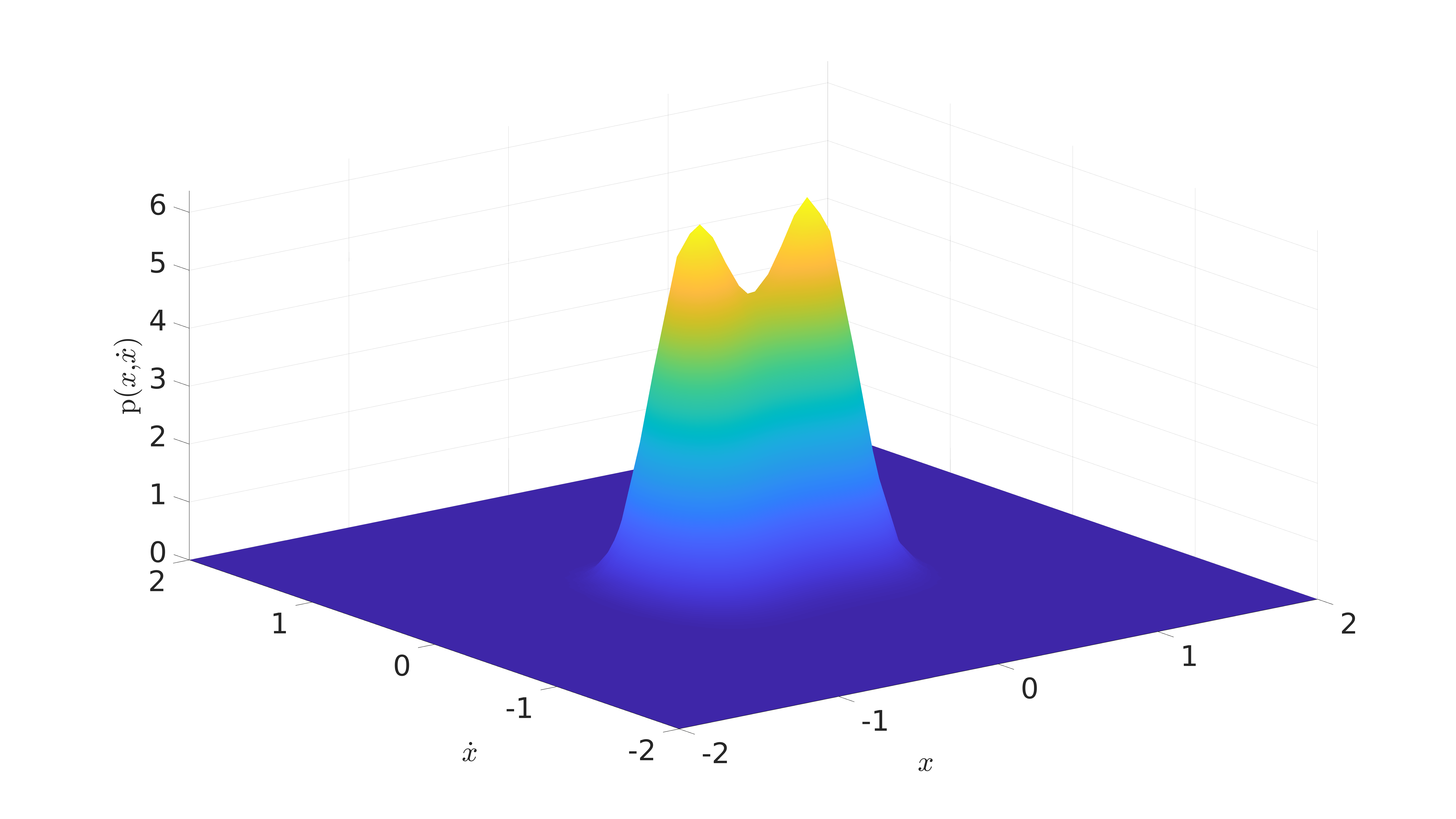}}
\subfigure[t = 10 s]{ \label{fig:DOcaseII_PDFc} \includegraphics[width=0.23\textwidth]{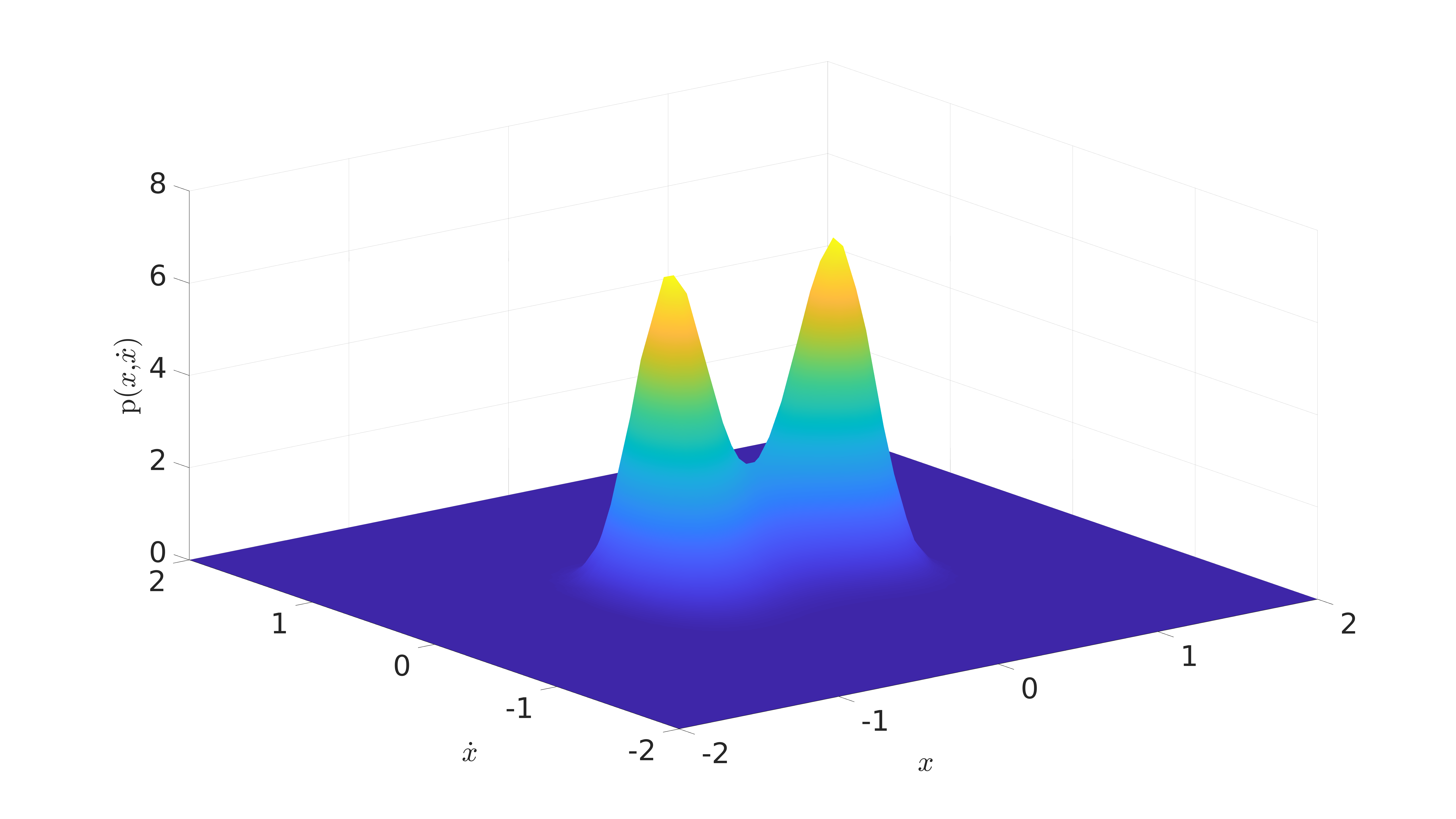}}
\subfigure[t = 50 s]{ \label{fig:DOcaseII_PDFd} \includegraphics[width=0.23\textwidth]{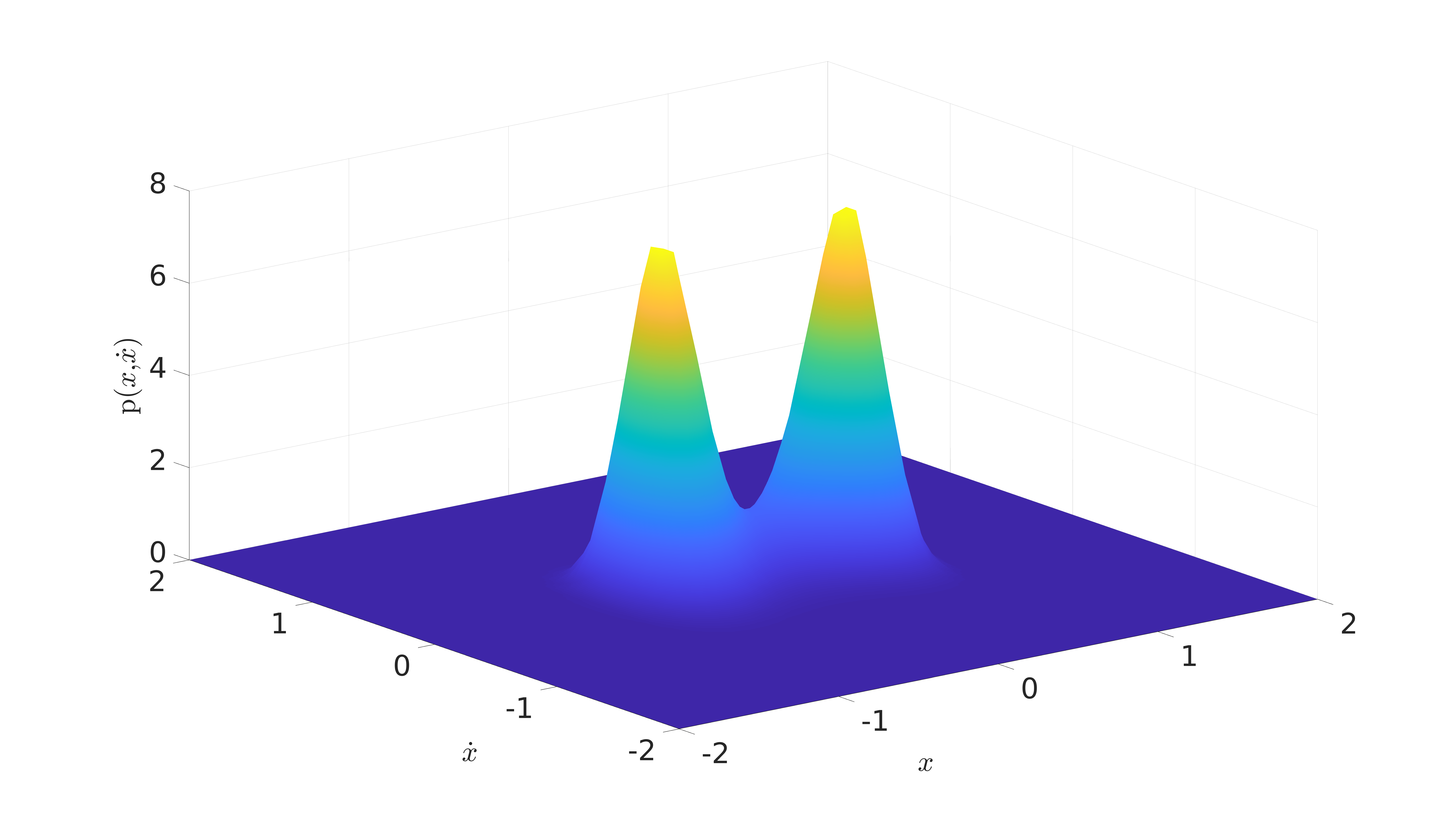}}
\caption{Duffing Oscillator : Pdf Surfaces}
\label{fig:DOcaseII_PDF}
\end{figure}

% \begin{figure}
% \centering
% \subfigure[Heat Map of Monomials]{\label{fig:DOcaseII_PDFa} \includegraphics[width=0.48\textwidth]{/src/Figures/DO_beta3/Heatmap_x.png}}
% \subfigure[Heat Map of Hamiltonians]{ \label{fig:DOcaseII_PDFb} \includegraphics[width=0.48\textwidth]{/src/Figures/DO_beta3/Heatmap_H.png}}
% \subfigure[Log-Pdf Ratio]{ \label{fig:DOcaseII_PDFc} \includegraphics[width=0.4\textwidth]{/src/Figures/DO_beta3/beta_ratio_H_x.png}}
% \caption{Duffing Oscillator : RS Coefficients Evolution and Log-Pdf Ratio.}
% \label{fig:DO_caseII_Heatmap_betaratio}
% \end{figure}
\begin{figure}
\centering
\subfigure[Heat Map of Monomials]{\label{fig:DOcaseII_PDFa} \includegraphics[width=0.48\textwidth]{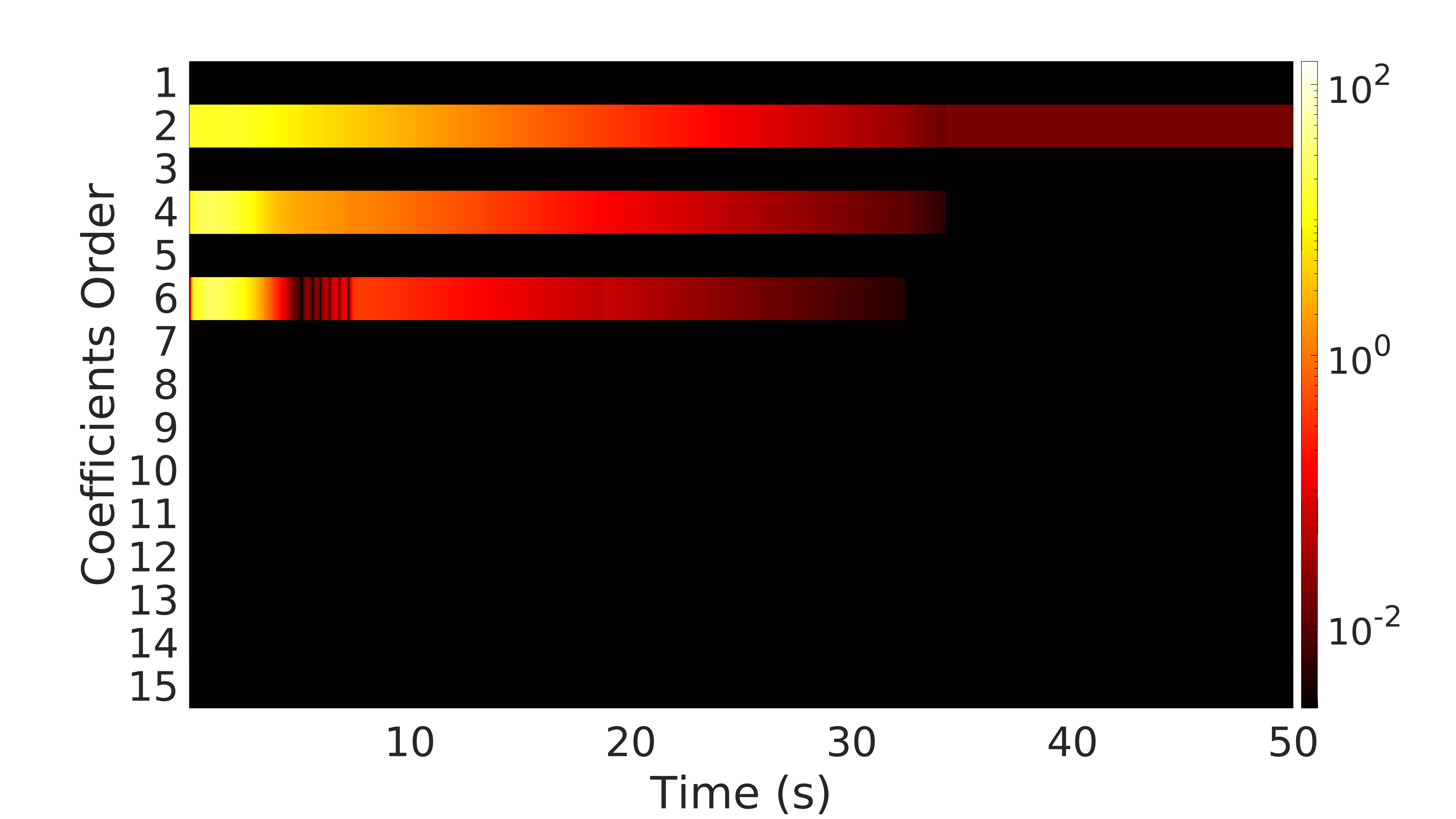}}
\subfigure[Heat Map of Hamiltonians]{ \label{fig:DOcaseII_PDFb} \includegraphics[width=0.48\textwidth]{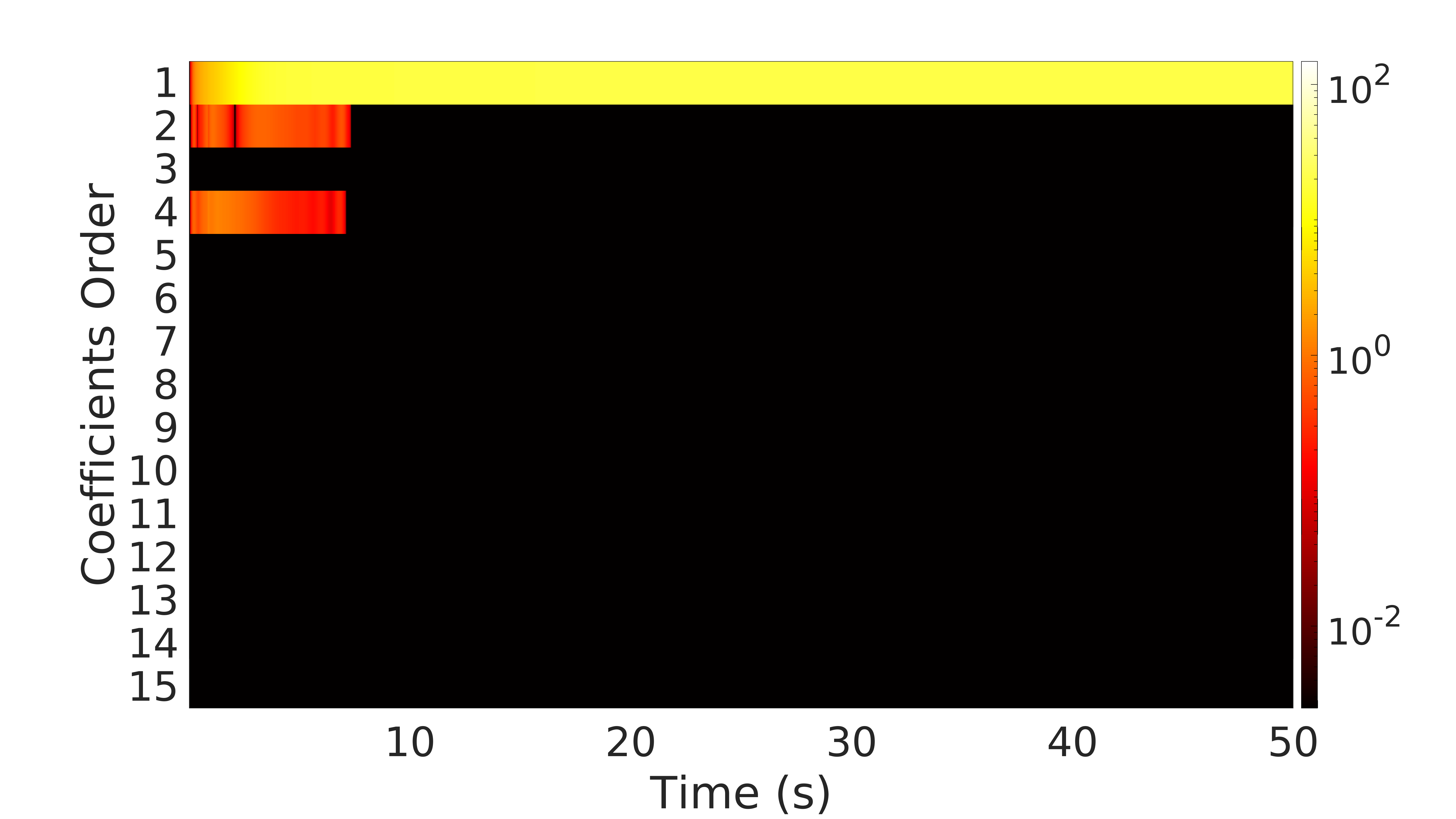}}
\caption{Duffing Oscillator : RS Coefficients Evolution}
\label{fig:DO_caseII_Heatmap}
\end{figure}

The first example corresponds to the noise-driven Duffing oscillator given in \eqref{DO_EOMs} and the stationary log-pdf is directly proportional to the Hamiltonian, H, and is shown in \eqref{Eq:DO_statPDF}. For simulation purposes, $Q = 1$, $\eta = 10$, $\alpha = -1$, $\beta = 3$ is considered. Furthermore, the Duffing oscillator with a soft spring ($\alpha \beta < 0$) and positive damping ($\eta >0$) has a bimodal stationary solution to FPKE. The modes of these bimodal solutions are concentrated at two stable equilibrium locations.

The initial pdf is assumed Gaussian, with zero mean and identity covariance. The basis function dictionary consists of monomials and Hamiltonians up to $15^{th}$ order resulting in total of $m = 151$ coefficients ($m_x$ = 136 and $m_h=15$). Further, 21 CUT8 collocation points are generated to sample the domain accurately. The RS parameters $\eta=1e-4$, $\Delta_s=1e-5$, and $\delta_{rs} = 1e-2$ is chosen. The time step $\Delta t = 0.01 s$ is considered and the system is solved for $t_f = 50 s$. The global domain is assumed to be $\x \in [-2,2]$, and the solution domain is mapped to a hypercube centered at the origin of length 2. In the local domain, a zero-mean Gaussian weighting function is employed with a covariance of $\boldsymbol{\Sigma}_w=\frac{1}{9} I_{2 \times 2}$ to ensure that the $\pm 3 \sigma$ bounds of the weight function are at the domain boundaries. 

The time-varying analysis is done by examining the pdf surfaces at different time steps. \reffig{fig:DOcaseII_PDF} shows the pdf surfaces of Duffing oscillator at time $t=2, 5, 10$ and $50s$. As the initial pdf was assumed to be Gaussian, it can be seen that the pdf at time $t=2s$ is near Gaussian, and as time progresses, the bimodal shape starts taking into effect from time $t=5s$. Finally, as time reaches $t=50s$, the true bimodal shape of Duffing oscillator is achieved. From this plot, the true stationary pdf is assumed to be given at $t=50s$. 

The convergence of the method can be studied further by examining the variation of the coefficients over time. \reffig{fig:DO_caseII_Heatmap} shows the heatmap of RS coefficients for the monomials and Hamiltonians, where the y-axis shows the same order coefficients are averaged at each time instance, and the x-axis shows the evolution with respect to time. It can be observed that the transient behavior of the pdf is dictated majorly by $2^{th}$, $4^{th}$, and $6^{th}$ order monomials, but as time progresses, these coefficients slowly start to diminish. Further, it can be noticed that as time increases, the Hamiltonian coefficient ($H^1$) slowly starts to become active and approaches the true stationary value of the order $1e2$. It is truly remarkable that the solution methodology automatically picks the actual Hamiltonian coefficient participating in the stationary solution out of 151 total basis functions.

\begin{figure}
% \hspace{1.25in}
\centering
\includegraphics[width=0.4\textwidth]{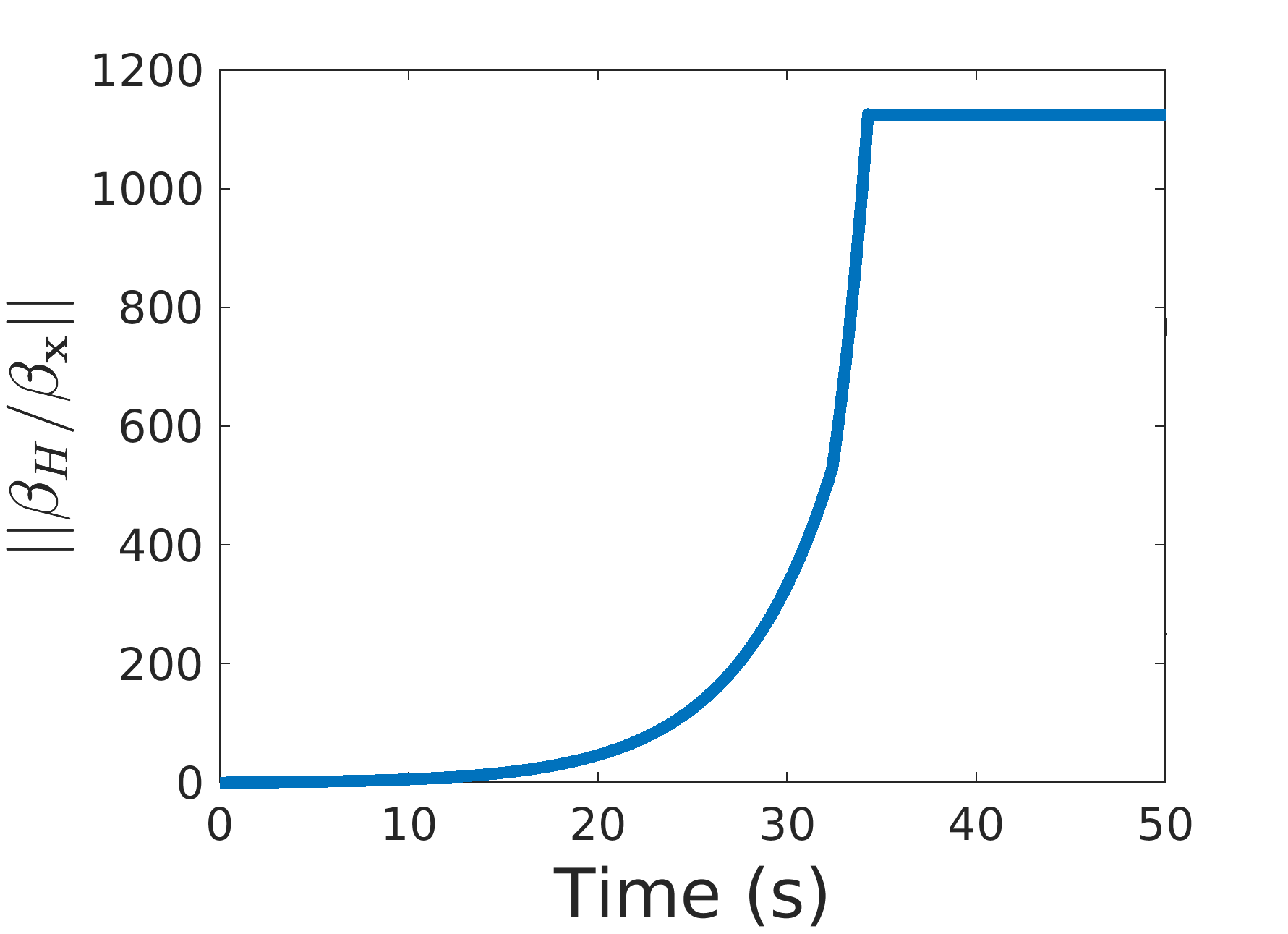}
\caption{Duffing Oscillator - Log-Pdf Ratio} \label{fig:DO_beta_ratio}
\end{figure}

Finally, in \reffig{fig:DO_beta_ratio}, the $\beta(H)/\beta(\x)$ which is the log-pdf ratio of Hamiltonians and monomials are shown. It is an important ratio to study as it provides the participation ratio of the Hamiltonian to the monomial basis functions in approximating the log-pdf. As time increases, the  $\beta(H)/\beta(\x)$ increases with time, signifying the increase in the dominance of Hamiltonian basis functions in computing the log-pdf. It can be further observed that the log-pdf ratio becomes constant after 35s, signifying that all the monomial coefficients become inactive.

\subsection{Two-Body Problem}

This section computes the uncertainty propagation for dynamical systems governed by the two-body equations of motion. Uncertainty in the system can be introduced due to errors in the initial condition or uncertain maneuvers. The sparse collocation method described above can be adapted for this type of non-linear system whose dynamics vary only as a function of the state. In this particular case, the evolution of the uncertainty can be found by propagating the initial state pdf through ordinary differential equations. As a result, an analytical solution can be found for the evolution of the state pdf. For an arbitrary orbital state-space, the system dynamics can be expressed as:
\begin{equation}\label{DynamicalSystem_Cont}
    \dot \x = \mathbf{f}(\x,t)
\end{equation} 
The state $\x_{k}$ at time $t_{k}$ can be written as a function of the initial state $\x_0$ using the system flow $\mathbf{F}$:
\begin{equation}\label{DynamicalSystem_Discrete}
    \mathbf{x}_{k} = \mathbf{F}( \x_0, t_{k}) 
\end{equation} 
Given the initial state pdf $p(\x_0)$, the state pdf $p(\mathbf{x}_k)$ at time $t_{k}$ can be computed using the transformation of variables as: 
\begin{equation}
p(\mathbf{x}_k) =p\left[\mathbf{x}_0=\mathbf{F}^{-1}(\mathbf{x}_k)\right]\left|\frac{\partial \mathbf{F}^{-1}}{\partial \mathbf{x}_k}\right|
\end{equation}
where flow $\mathbf{F}$ is assumed to be invertible, continuously differentiable mapping, with inverse given by: $\x_0=\mathbf{F}^{-1}(\mathbf{x}_k)$. This enables the determination of the true state pdf from the initial state pdf and the system flow information. If the dynamics are precisely known, it is possible to simplify further the above equation, which can be written as:
\begin{equation}\label{Eq:statepdf_xdot_fx}
p(\mathbf{x}_k)=p\left[\mathbf{x}_0=\mathbf{F}^{-1}(\mathbf{x}_k)\right] \exp \left(-\int_{0}^{t_k} \nabla \cdot \mathbf{f}(\mathbf{x}(\tau), \tau) d \tau\right)
\end{equation}
Using the Method of Characteristics, the aforementioned result can be obtained. This equation is significant because it allows one to determine the discrete-time probability value $p(\mathbf{x}_k)$ of any initial pdf sample. A particular case in propagating the uncertainty of a two-body problem due to maneuver uncertainties can also be considered similarly. For the computation of a numerical solution, the initial collocation points corresponding to the post-maneuver mean $\mu_{\x_1}$ and covariance $\Sigma_{\x_1}$ must be propagated through the dynamics, and the corresponding state pdf needs to be calculated. It should be noted that the propagation of initial collocation points automatically tracks the solution domain of the FPKE. 

A solution approach using \eqref{Eq:statepdf_xdot_fx} is sought to compute the numerical approximation of this system. Due to the high-dimensional input space and to enforce the smoothness constraint, the state pdf is assumed to have the following structure:

\begin{equation}
p\left(\mathbf{x}_{k+1}, t_{k+1}\right)=\delta p\left(\mathbf{x}_{k+1}, t_{k+1}\right) p\left(\mathbf{x}_{k+1}, t_{k}\right)
\end{equation}
where $p\left(\mathbf{x}_{k+1}, t_{k+1}\right)$ is the pdf of state $\x_{k+1}$ at time $t_{k+1}$, $p\left(\mathbf{x}_{k+1}, t_{k}\right)$ is the pdf of state $\x_{k+1}$ at time $t_{k}$, and $\delta p\left(\mathbf{x}_{k+1}, t_{k+1}\right)$ is the departure pdf of state $\x_{k+1}$ from time $t_{k}$ to $t_{k+1}$. Here the state pdf at previous time step is used as the weight function, i.e., $p_W=p\left(\mathbf{x}_{k+1}, t_{k}\right)$ to impose the infinite boundary conditions. Transforming the above equation into log-pdf and assuming a polynomial basis function results in:

\begin{equation}
\begin{aligned}
\beta\left(\mathbf{x}_{k+1}, t_{k+1}\right)= &\delta \beta\left(\mathbf{x}_{k+1}, t_{k+1}\right)  + \beta\left(\mathbf{x}_{k+1}, t_{k}\right) \\
\ln \left[p\left(\mathbf{x}_{k+1}, t_{k+1}\right)\right] \approx \mathbf{c}_{k+1}^{T} \Phi\left(\mathbf{x}_{k+1} \right)= &\delta \mathbf{c}_{k+1}^{T} \Phi\left(\mathbf{x}_{k+1}\right)+\mathbf{c}_{k}^{T} \Phi\left(\mathbf{x}_{k+1}\right) \\
\end{aligned}\label{eq:del_ck}
\end{equation}
where $\mathbf{c}_{k+1}$ is a vector of unknown coefficients at time $t_{k+1}$, $\delta\mathbf{c}_{k+1}$ is a vector of unknown departure coefficients from time $t_{k}$ to time $t_{k+1}$, $\mathbf{c}_{k}$ is the coefficients at time $t_{k}$ and $\Phi(\x_{k+1})$ is a vector of the chosen basis functions of state $\x_{k+1}$. The objective is first to compute the departure coefficients $\delta\mathbf{c}_{k+1}$ that can serve as a weight function for the previous time pdf and thus enforce the smoothness constraint. Using the aforementioned equation, the coefficients $\mathbf{c}_{k+1}$ can be found as:
\begin{equation}
\mathbf{c}_{k+1}  =  \delta \mathbf{c}_{k+1}   +\mathbf{c}_{k}  
\end{equation}
To find the optimal coefficients $\delta\mathbf{c}_{k+1}$, Eq. \ref{eq:del_ck} can be rearranged as:
\begin{equation}
\delta \mathbf{c}_{k+1}^{T} \Phi\left(\mathbf{x}^j_{k+1}\right)=\ln \left[p\left(\mathbf{x}^j_{k+1}\right)\right] - \mathbf{c}_{k}^{T} \Phi\left(\mathbf{x}^j_{k+1}\right)
\end{equation}
Further simplifying the above equation results in:
\begin{equation}
\mathbf{A} \delta \mathbf{c}_{k+1}=\mathbf{b} 
\end{equation}
where:
% \begin{equation}
% \begin{gathered}
% \mathbf{A}_{j}=\Phi^{T}\left(\mathbf{x}^j_{k+1}\right); \qquad \mathbf{b}_{j} =\ln \left[p\left(\mathbf{x}^j_{k+1}\right)\right]-\Phi^{T}\left(\mathbf{x}^j_{k+1}\right)  \mathbf{c}_{k}, \quad j=1,2, \ldots, N
% \end{gathered}
% \end{equation}
\begin{equation}
\begin{gathered}
\mathbf{A}_{j}=\Phi^{T}\left(\mathbf{x}^j_{k+1}\right), \quad j=1,2, \ldots, N \\
\mathbf{b}_{j} =\ln \left[p\left(\mathbf{x}^j_{k+1}\right)\right]-\Phi^{T}\left(\mathbf{x}^j_{k+1}\right)  \mathbf{c}_{k}, \quad j=1,2, \ldots, N
\end{gathered}
\end{equation}
Therefore, the following optimization problem is presented to solve for the departure coefficients $\delta \mathbf{c}_{k+1}$:
$$
\min _{\delta \mathbf{c}_{k+1}}\left\|\mathbf{K} \delta \mathbf{c}_{k+1}\right\|_{1}
$$
subject to: $$
\left\|\mathbf{W}(\mathbf{A} \delta \mathbf{c}_{k+1}-\mathbf{b}) \right\|_{2} \leq \boldsymbol{\epsilon}
$$
Note that the departure from the previous time pdf must be solved by calculating the departure coefficients, $\delta \mathbf{c}_{k+1}$. The procedure for obtaining the $l_2$-norm, $l_1$-norm, and RS solution for determining the coefficients remain the same as illustrated in Algorithm \ref{SparseSolution_algo}.

\begin{figure}
\centering
\subfigure[Case - I : Monomials]{\label{fig:TB_coeff_Heatmapa} \includegraphics[width=0.32\textwidth]{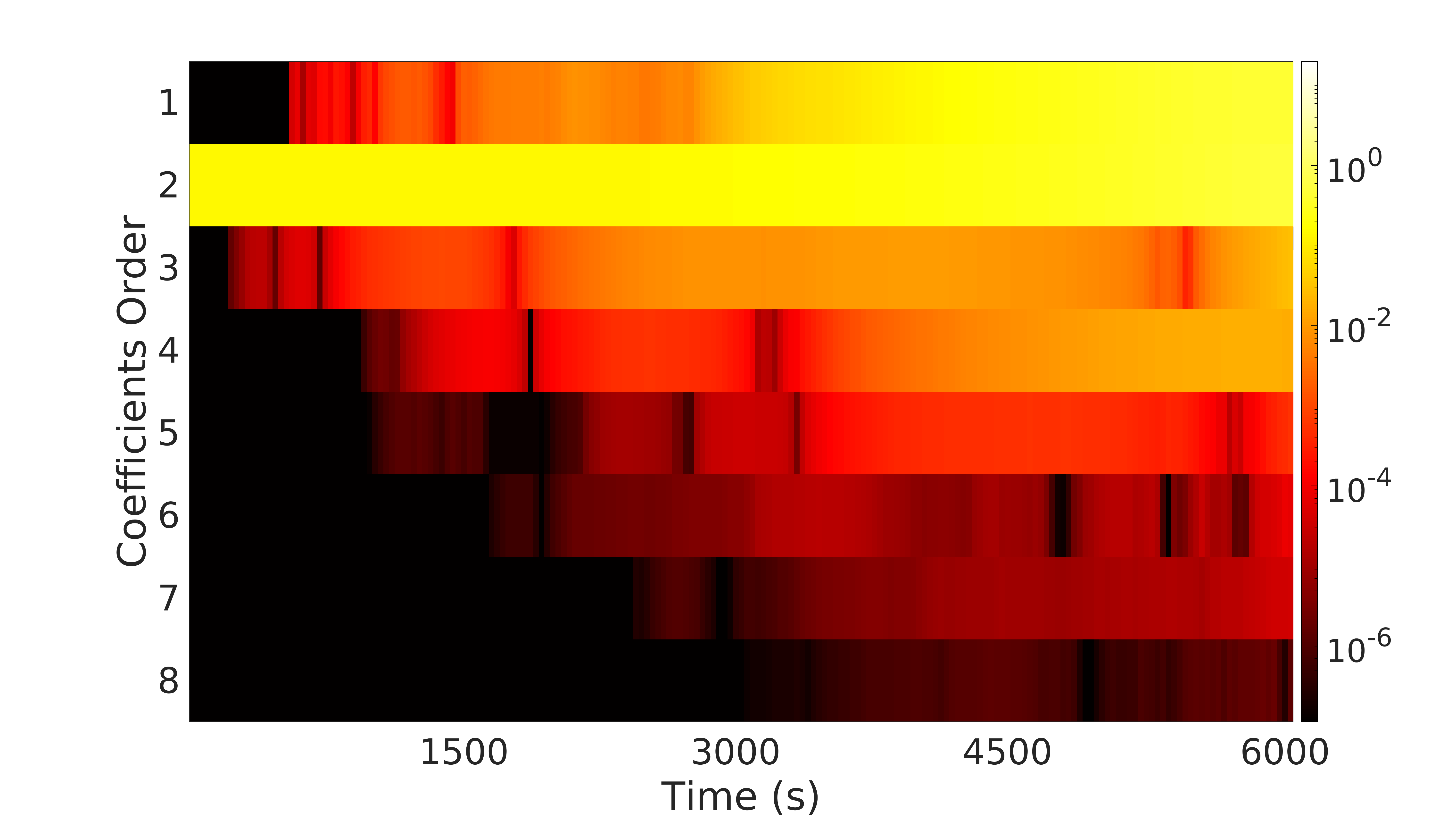}}
\subfigure[Case - II : Monomials]{ \label{fig:TB_coeff_Heatmapb} \includegraphics[width=0.32\textwidth]{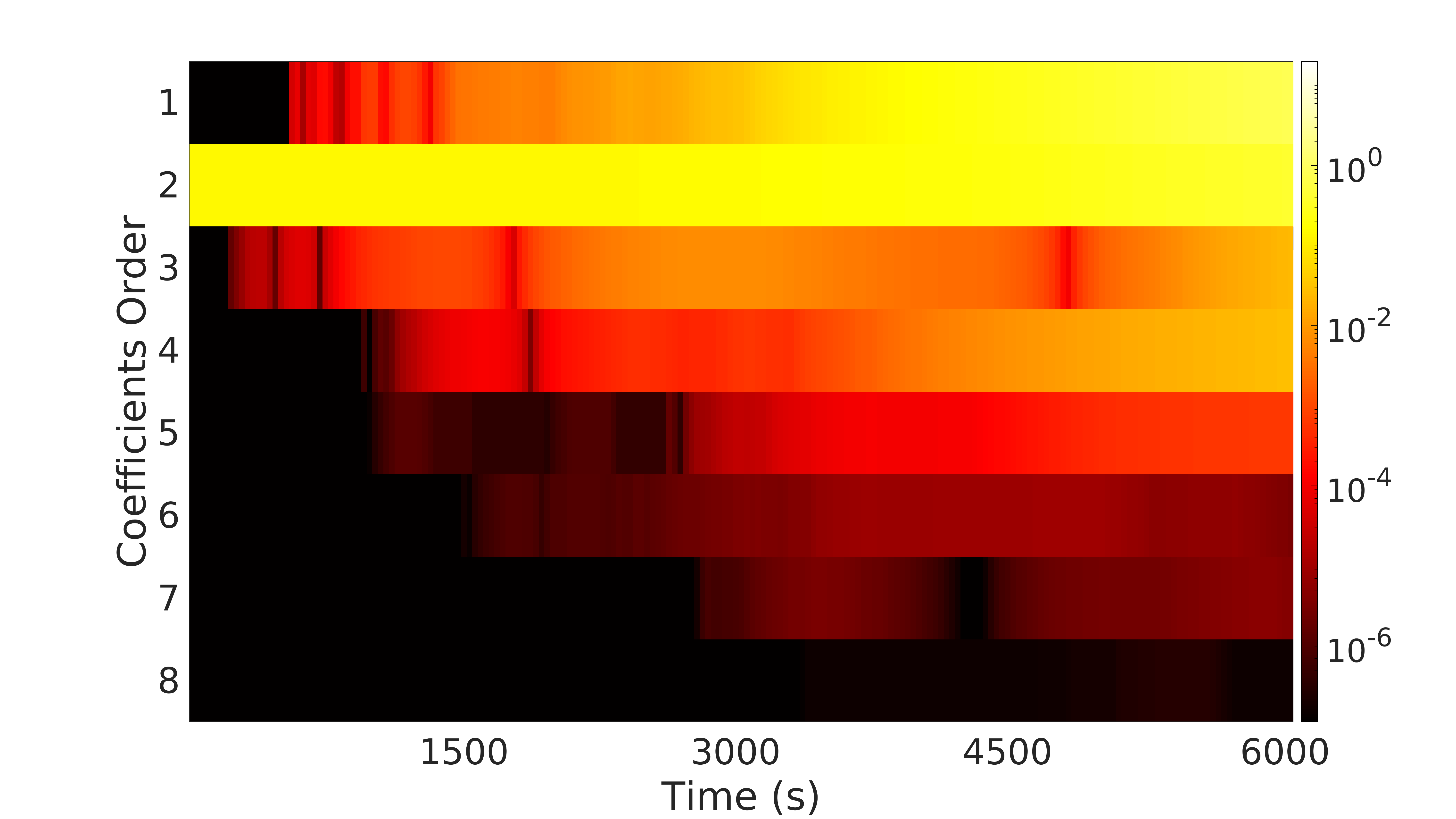}}
\subfigure[Case - II : Hamiltonians]{ \label{fig:TB_coeff_Heatmapc} \includegraphics[width=0.32\textwidth]{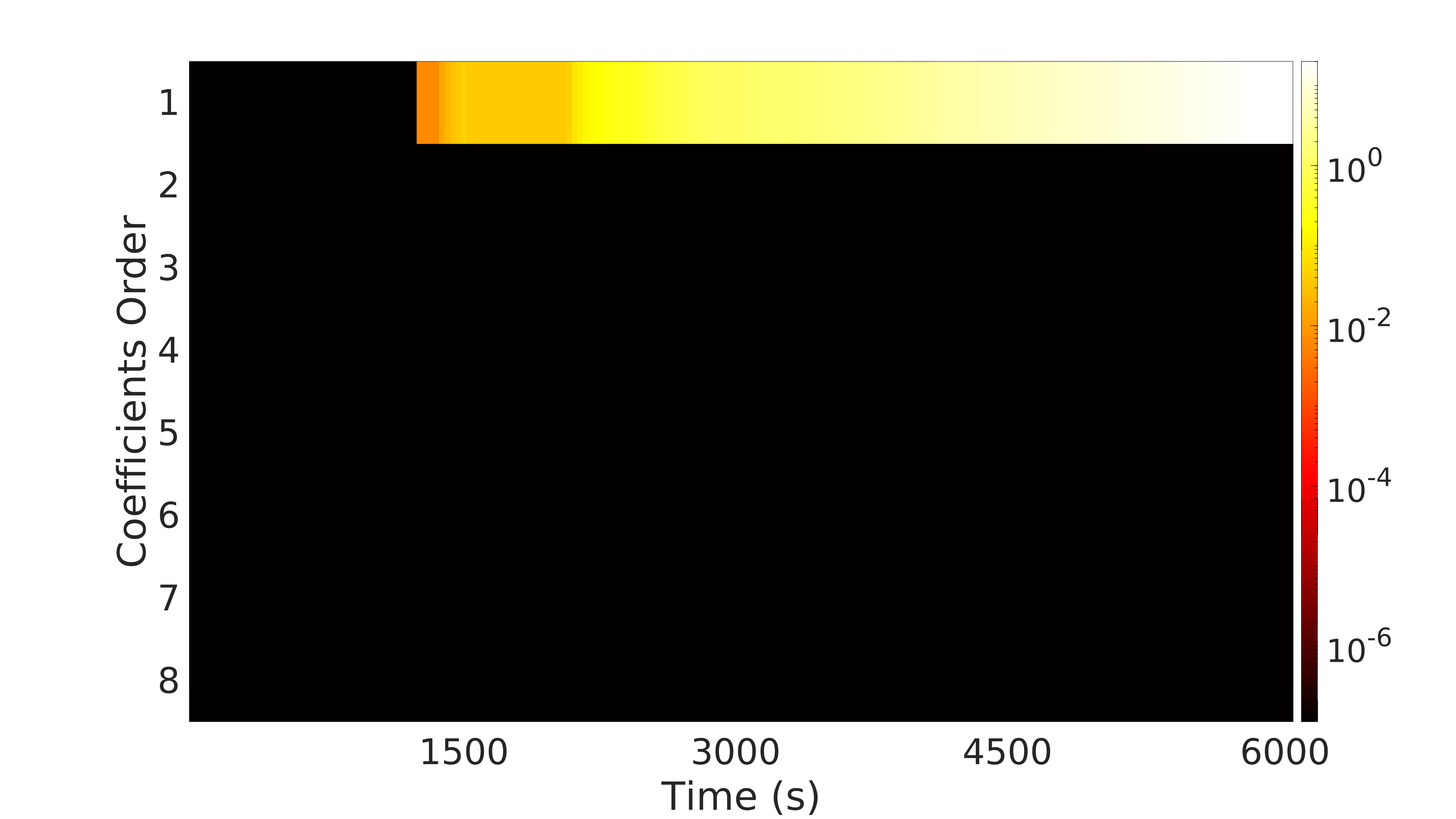}}
\caption{Orbit Transfer Maneuver - Heat Map of RS Coefficients}
\label{fig:TB_coeff_Heatmap}
\end{figure}

% \begin{figure}
% \centering
% \subfigure[Case - I : Monomials]{\label{fig:TB_coeff_Heatmapa} \includegraphics[width=0.32\textwidth]{/src/Figures/TB_x/Heatmap_log.png}}
% \subfigure[Case - II : Monomials]{ \label{fig:TB_coeff_Heatmapb} \includegraphics[width=0.32\textwidth]{/src/Figures/TB_H_x/Heatmap_log_x.png}}
% \subfigure[Case - II : Hamiltonians]{ \label{fig:TB_coeff_Heatmapc} \includegraphics[width=0.32\textwidth]{/src/Figures/TB_H_x/Heatmap_log_H.png}}
% \caption{Orbit Transfer Maneuver - Heat Map of RS Coefficients}
% \label{fig:TB_coeff_Heatmap}
% \end{figure}

\subsubsection{Orbit Transfer Maneuver}

This section describes the numerical analysis of a cooperative maneuvering satellite used in ~\cite{betts1977optimal,hall2019probabilistic}. This satellite illustrates the optimal three-burn transfer orbit from Low Earth Orbit (LEO) to a final orbit with a particular inclination. For the three burns, it is assumed that the nominal target maneuver of the satellite is known. For the first burn, the uncertainty is propagated by assuming uncertain initial conditions and an uncertain burn profile. This procedure can be repeated to propagate the uncertainty for the subsequent burns. In addition, it is assumed that the system's uncertainty is Gaussian, with the mean given by the nominal profile. The equation of motion for this example is given as:
\begin{equation}
\begin{aligned}
\dot \x = & f(\x, t) + g(\mathbf{u}, t) =  \left[ \begin{array}{c}
\mathbf{v}\\
-\frac{\mu \mathbf{r}}{|\mathbf{r}|^{3}}
\end{array} \right] + \left[ \begin{array}{c}
0 \\
\mathbf{T} \mathbf{u} 
\end{array} \right] \\
\end{aligned}
\end{equation}
where $f(\x, t)$ is given by two-body equations of motion and $g(\mathbf{u},t)$ is a control maneuver model that is only active when a maneuver must be performed. A particular case for this system can be considered when there is no maneuver required to perform, i.e., $g(\mathbf{u},t) =0$. The system dynamics is then reduced to the two-body equation of motion. The control input for this impulsive maneuver can be written in the Earth-Centered Inertial (ECI) frame as follows:
\begin{equation}
\mathbf{u}=\left[\begin{array}{c}
|\Delta \mathbf{v}| \sin \theta \cos \psi \hat{N} \\
|\Delta \mathbf{v}| \cos \theta \cos \psi \hat{T} \\
-|\Delta \mathbf{v}| \sin \psi \hat{W}
\end{array}\right]
\end{equation}
where $|\Delta v|$ represents the magnitude of velocity, $\psi$ represents the pitch, and $\theta$ represents the yaw angles in the NTW satellite coordinate system. 

\begin{table}[htb!]
\vspace{0.1in}
    \begin{minipage}{.5\linewidth}
      \centering
        \begin{tabular}{|c|c|c|c|}
            \hline Parameter & Parking Orbit & $1^{\text {st }}$ Transfer \\
            \hline$a^{*}(k m)$ & $6667.32$ & $22,835.4$   \\
            \hline$e^{*}$ & 0 & $0.7080$   \\
            \hline$i^{*}\left(^{\circ}\right)$ & $37.40$ & $35.78$   \\
            \hline$\omega^{*}\left(^{\circ}\right)$ & N/A & $254.0$   \\
            \hline
        \end{tabular}
        \vspace{0.15in}
      \caption*{Transfer Orbit Parameters\label{table:OrbitParameter}}
    \end{minipage}%
    \begin{minipage}{.5\linewidth}
      \centering
      \vspace{-0.18in}
        \begin{tabular}{|c|c|c|c|}
            \hline Parameter & $1^{\text {st }}$ Burn & Final Position   \\
            \hline$\left|\Delta v^{*}\right|(m / s)$ & $2,383.5$ &$\mathrm{~N} / \mathrm{A}$ \\
            \hline$\theta^{*}\left(^{\circ}\right)$ & $-1.36$ &   $\mathrm{~N} / \mathrm{A}$ \\
            \hline$\psi^{*}\left(^{\circ}\right)$ & $5.99$ &  $\mathrm{~N} / \mathrm{A}$ \\
            \hline$\eta^{*}\left(^{\circ}\right)$ & 255 &   $121.7$ \\
            \hline$t(s)$ & 0 & $6,078.2$   \\
            \hline
        \end{tabular}
        \vspace{0.10in}
        \caption*{Burn Parameters\label{table:BurnParameter}}
    \end{minipage} 
    \caption{Nominal Transfer Orbit and Burn Parameters\label{table:NominalParameter}}
\end{table}

The nominal transfer orbital elements and burn parameters used in the circular parking orbit are listed in Table \ref{table:NominalParameter}. The first burn is applied to the parking orbit in LEO with $i=37.4^{\circ}$ at a latitude argument of $u=255^{\circ}$. The nominal orbital elements for the transfer orbits are listed on the left side of Table \ref{table:NominalParameter}, with the first burn occurring at $t_1=0$ and the second at $t_2=6078.2 s$. The mean of the initial state corresponds to the parking orbit with the argument of latitude $\eta=255^{\circ}$. The standard deviation of $50 m$ is applied to the position, and $0.1 m/s$ is applied to the velocity in each direction. The nominal control parameters can be calculated based on the nominal burn parameters on the right side of Table \ref{table:NominalParameter}. The standard deviation of $5 m/s$ is applied to the impulsive control velocity in each direction.

For this problem, uncertainty propagation can be separated into two parts: i) state uncertainty propagation in the presence of maneuver uncertainty $g(\mathbf{u},t)$ and ii) state uncertainty propagation when $g(\mathbf{u},t)=0$. In the first part, the uncertain instantaneous maneuver is applied to the uncertain initial conditions. Since the burn is instantaneous, the satellite's position remains constant, and the uncertainty propagates only in the velocity vector. Also, utilizing the fact that the initial velocity and applied control input are considered Gaussian, the parameters of propagated uncertainty in the velocity can be computed using Gaussian properties. Let $\mathbf{v}_{0}$ represent the initial velocity, $\mathbf{u}_{0}$ represent the control input, and $\mathbf{v}_{1}$ represent the propagated velocity. The mean $\mu_{\mathbf{v}_1}$ and covariance $\Sigma_{\mathbf{v}_1}$ of the propagated velocity can therefore be computed as follows:
\begin{equation}
\begin{aligned}
\mathbf{v}_1=& \mathbf{v}_{0} + \mathbf{T} \mathbf{u}_0 \\
\mu_{\mathbf{v}_1}= & E[ \mathbf{v}_1]  = E[\mathbf{v}_{0}] + \mathbf{T} E[\mathbf{u}_0] \\
\Sigma_{\mathbf{v}_1} = & E[\mathbf{v}_1\mathbf{v}_1^T] =  E[(\mathbf{v}_{0} + \mathbf{T} \mathbf{u}_0)(\mathbf{v}_{0} + \mathbf{T} \mathbf{u}_0)^T] =  E[\mathbf{v}_{0} \mathbf{v}_{0}^T] + \mathbf{T} E[\mathbf{u}_0 \mathbf{u}^T_0]\mathbf{T}^T \\
\end{aligned}
\end{equation}
Finally, the uncertain parameters of the post-maneuver state $\x_1$ can be calculated by concatenating the mean and covariance of the initial position $\x_0$ and the propagated velocity, $\mathbf{v}_1$. The post-maneuver mean and covariance of $\x_1$ can be used to find the state pdf $p(\x_1)$. The second component is the propagation of this uncertainty through the two-body equations of motion, which were previously derived.

For computing the numerical solution of this example, two different test cases are considered to see the effect of including Hamiltonian coefficients in approximating the state pdf. In test case - I, the monomial basis functions up to $8^{th}$ order are considered in the dictionary, which results in $m =3003$ coefficients. Whereas in test case - II, monomials and Hamiltonians up to $8^{th}$ order are considered, resulting in total $m =3011$ coefficients. For comparison, all other parameters are kept constant for both test cases. The CUT8 points are chosen as collocation points for this $6-D$ system, resulting in 745 points. The parameters $\alpha = 1e-6$, $\eta=1e-4$, $\Delta_s=1e-4$, and $\delta_{rs} = 1e-4$ is chosen.

\begin{figure}
\centering
\subfigure[Case - I ]{\label{fig:Evol_NNZa} \includegraphics[width=0.45\textwidth]{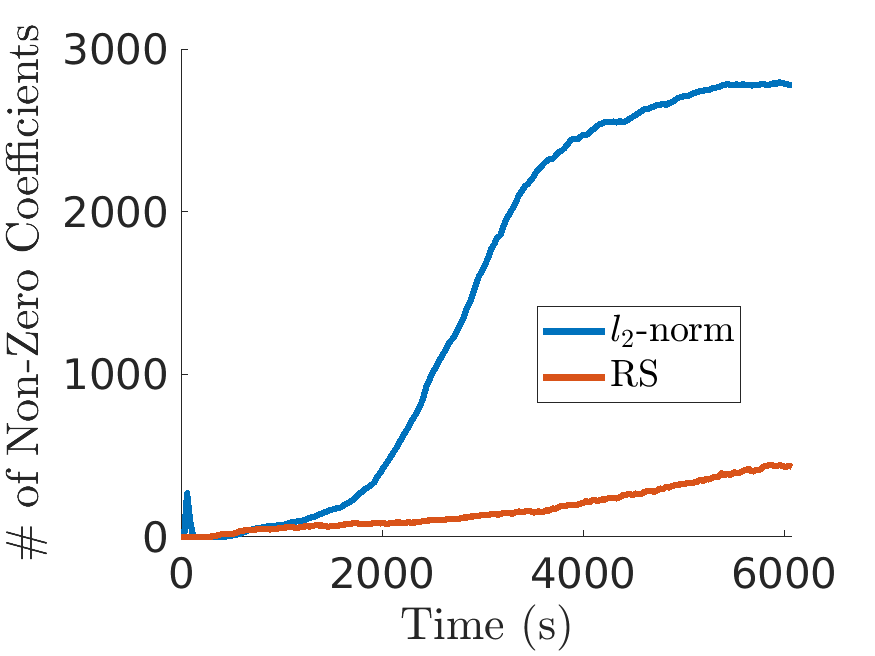}}
\subfigure[Case - II ]{ \label{fig:Evol_NNZb} \includegraphics[width=0.45\textwidth]{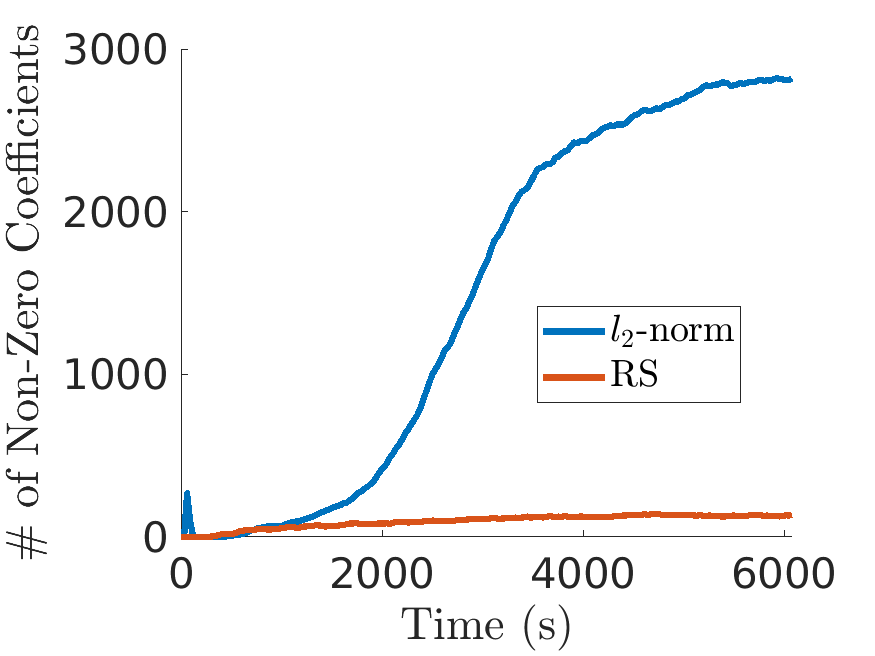}}
\caption{Orbit Transfer Maneuver - Evolution of Non-Zero Coefficients}
\label{fig:Evol_NNZ}
\end{figure}

\reffig{fig:TB_coeff_Heatmap} illustrates the heatmap of the RS coefficients (averaged over the similar order) for the test case - I and test case - II in the left and right columns, respectively. By comparing the monomials color bar in both cases, it can be seen that coefficients up to $5^{th}$ order are active in both cases. However, the maximum value of the higher order coefficients from $6^{th}$ to $8^{th}$ order in case - I comes out to be of the order $1e-1$ while in case - II, it comes out to be of the order $1e-3$. This is because the first order Hamiltonian coefficient ($H^1$) in case - II replaces most of the monomial coefficients and can quantify the uncertainty with fewer coefficients. The heatmap of Hamiltonians provides the information that the Hamiltonian of order one ($H^1$) is only active, and the rest Hamiltonians are inactive. Thus, including the Hamiltonian dictionary of basis function performs better and requires less number of coefficients for uncertainty propagation.

% the maximum coefficient magnitude in case - II is less than in case - I. Moreover, by looking carefully, it can be observed that in case - I, the monomials up to $5^{th}$ order are active, whereas, in case - II, the monomials up to $4^{th}$ order are active. This is because the first order Hamiltonian coefficient ($H^1$) in case - II is replacing most of the monomial coefficients and can approximate the reachability set with fewer coefficients. The heatmap of Hamiltonians provides the information that the Hamiltonian of order one ($H^1$) is only active, and the rest Hamiltonians are inactive. \reffig{fig:TB_coeff678} compares the magnitude of the higher order monomials for both cases. Clearly, it can be seen that by including the Hamiltonian in the dictionary of basis function, the magnitude of higher order basis functions are of the order $1e-3$, whereas in case-I, it comes out to be of the order $1e-1$. Thus, including the Hamiltonian dictionary of basis function performs better and requires less number of coefficients in computing the reachability analysis.

\reffig{fig:Evol_NNZ} shows the enforced sparsity constraint for both cases where the evolution of the number of non-zero coefficients for the $l_2$-norm and RS coefficients is shown. Non-zero coefficients are defined here as those whose magnitude is greater than $1e-5$. As mentioned earlier, 3003 and 3011 are the total number of coefficients for case - I and case - II, respectively. It can be noticed here that the non-zero coefficients required by the $l_2$-norm solution increase exponentially with time before converging to approximately 2600 coefficients. In contrast, the RS solution requires around 500 coefficients for case - I and less than 100 coefficients for case - II at any given time. Again, by including the Hamiltonian $H^1$ in the dictionary of basis function, about 400 fewer coefficients are needed.

% \begin{figure}
%  \centering
%  \begin{tabular}{@{}p{2in}  @{}p{2in}} 
%   \addheight{\includegraphics[width=1.8in]{/src/Figures/TB_x/NormErr_pdf_Normalzed.png}} &  \hspace{0.15in} {\includegraphics[width=1.8in]{/src/Figures/TB_H_x/NormErr_pdf_Normalzed.png}} \\ 
%   \addheight{\includegraphics[width=1.8in]{/src/Figures/TB_x/TwoNormErr_Normalzed_tvec.png}} & \hspace{0.15in} {\includegraphics[width=1.8in]{/src/Figures/TB_H_x/TwoNormErr_Normalzed_tvec.png}}  \\ 
%   \vspace{0.5in} \\
%   \hspace{0.8in} (a) Case - I & \hspace{0.8in} (b)  Case - II \\
% \end{tabular}
%  \caption{Orbit Transfer Maneuver - Training vs. Testing Pdf and Log-Pdf Error using RS Coefficients} \label{fig:TestTrain_caseI_II}
%  \end{figure}

% \reffig{fig:TestTrain_caseI_II} depicts the normalized two-norm error of log-pdf and state pdf for both cases with training and testing data. $1,000$ MC samples are taken at each time step for testing data. In each of these cases, the training error comes out to be less than the testing error. However, for case - I, it can be noticed that the state pdf error for the testing data is going to enormous values in the last few timesteps while the log-pdf error is showing normal behavior. This numerical anomaly is likely due to the involvement of a large number of coefficients in approximating the pdf, which again shows that including Hamiltonians in the dictionary of basis functions avoids these numerical issues and, in turn, improves the accuracy of the reachability analysis. 

\begin{figure} 
 \centering
\includegraphics[width=0.35\textwidth]{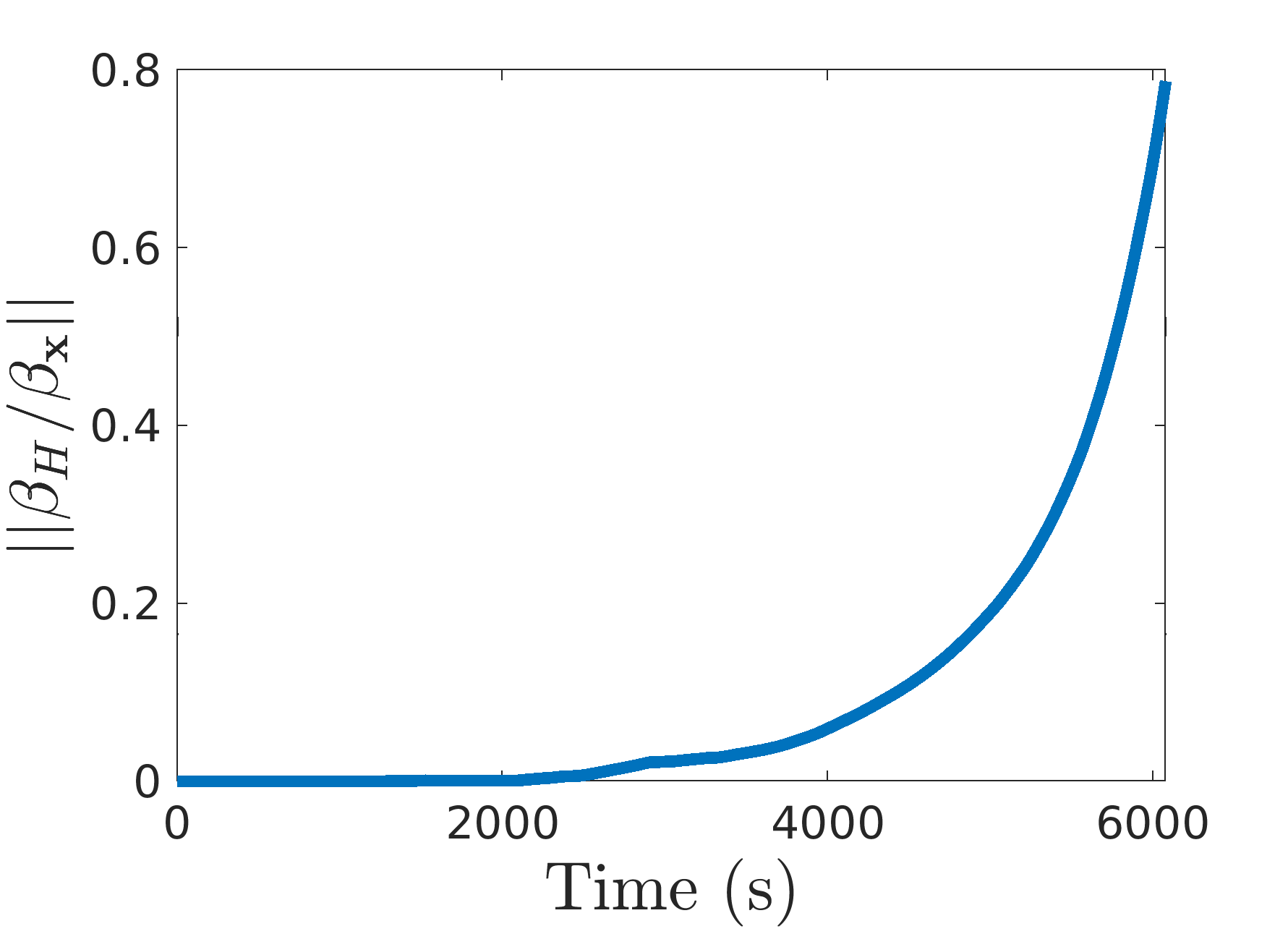}
\caption{Orbit Transfer Maneuver Case - II - Log-Pdf Ratio} \label{fig:beta_ratio}
\end{figure}

\begin{figure}
\centering
\begin{tabular}{@{}M{0.3in}@{}p{\spcing}@{}p{\spcing}@{}p{\spcing}@{}p{\spcing}} 
{\multirow{-5.5}{*}{\rotatebox[origin=c]{90}{\scriptsize  Time = $t_f/2$  }}} 
& \addheight{\includegraphics[width=1.3in]{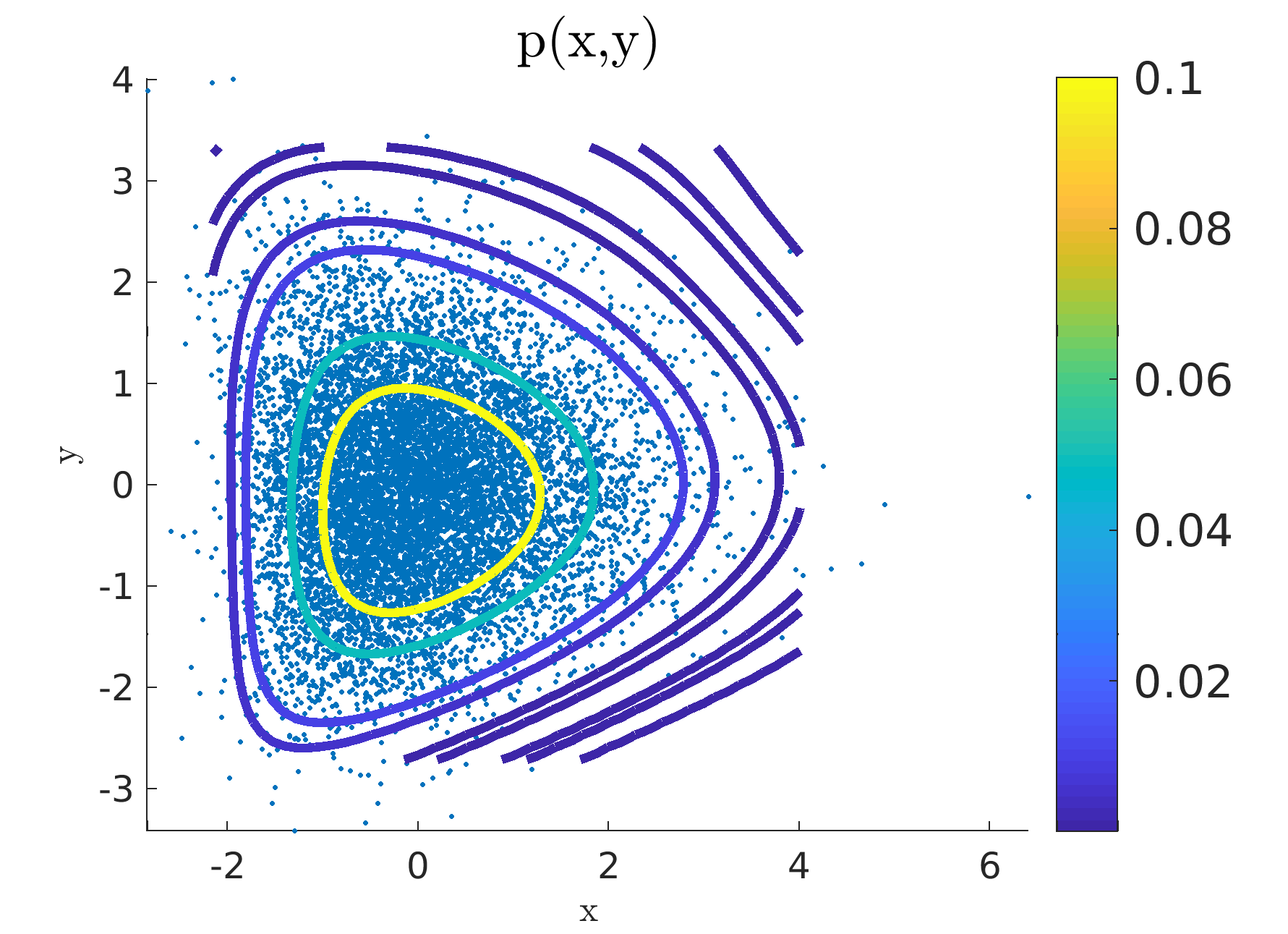}} 
& {\includegraphics[width=1.3in]{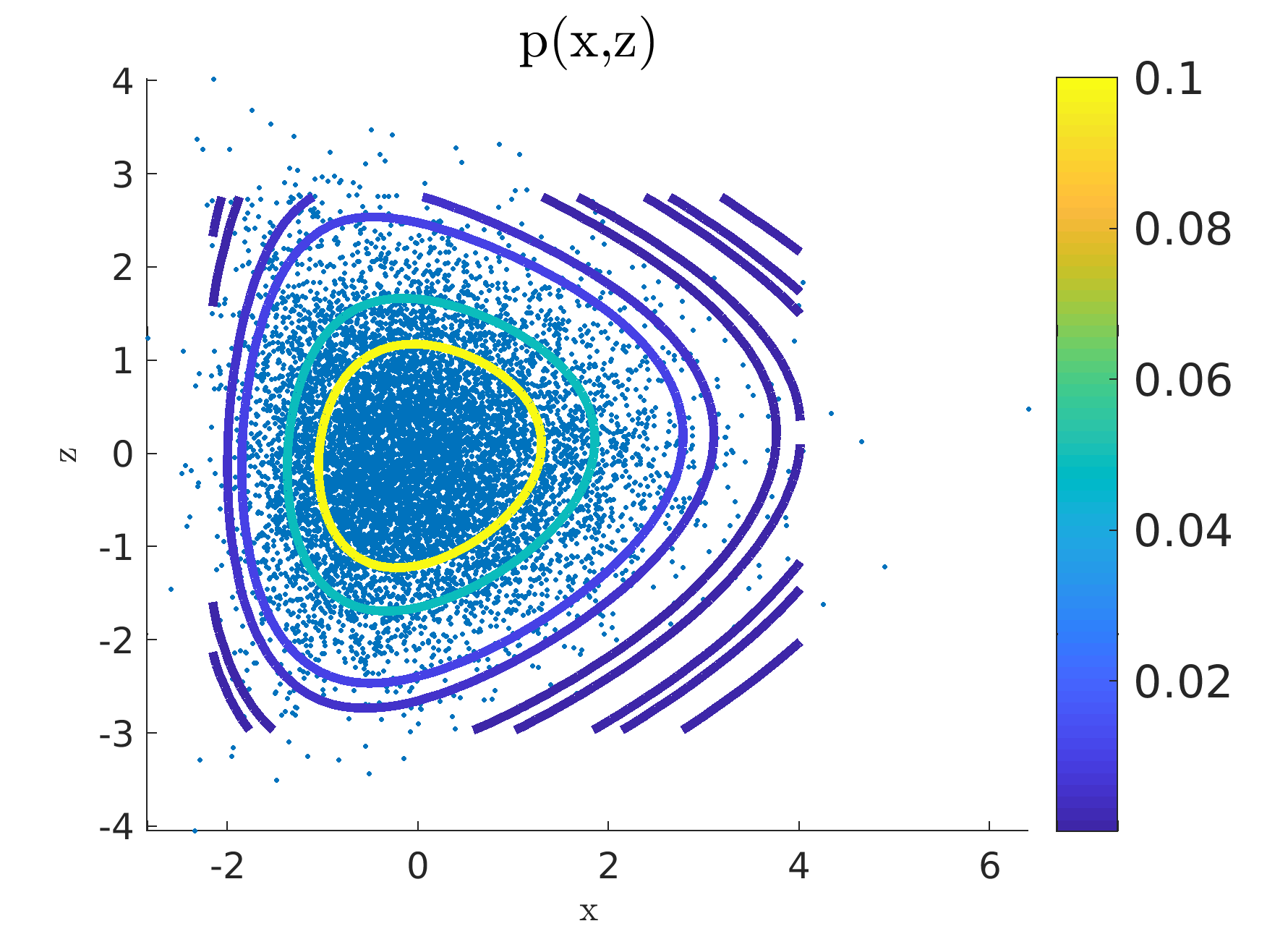}} 
& {\includegraphics[width=1.3in]{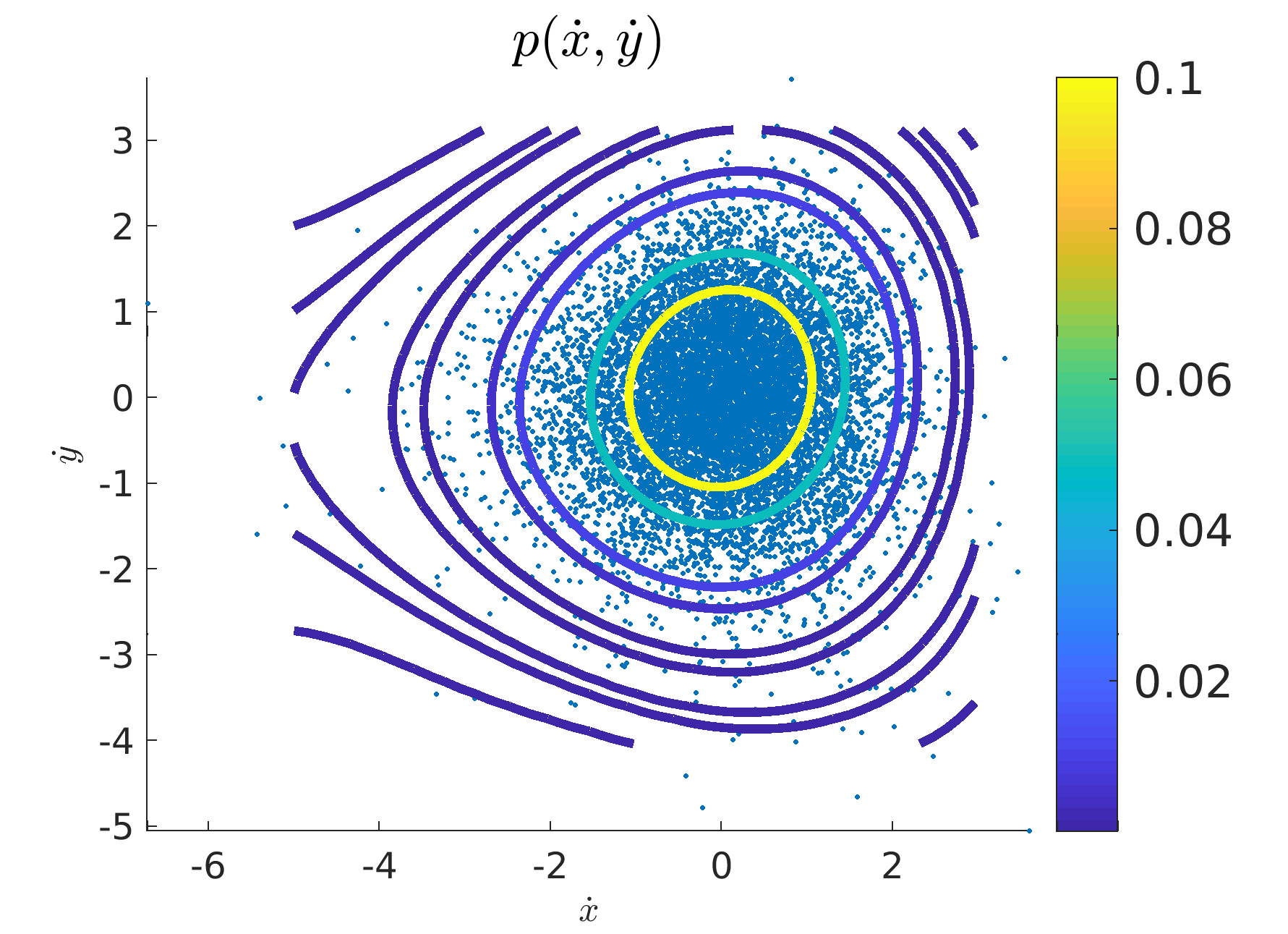}} 
& {\includegraphics[width=1.3in]{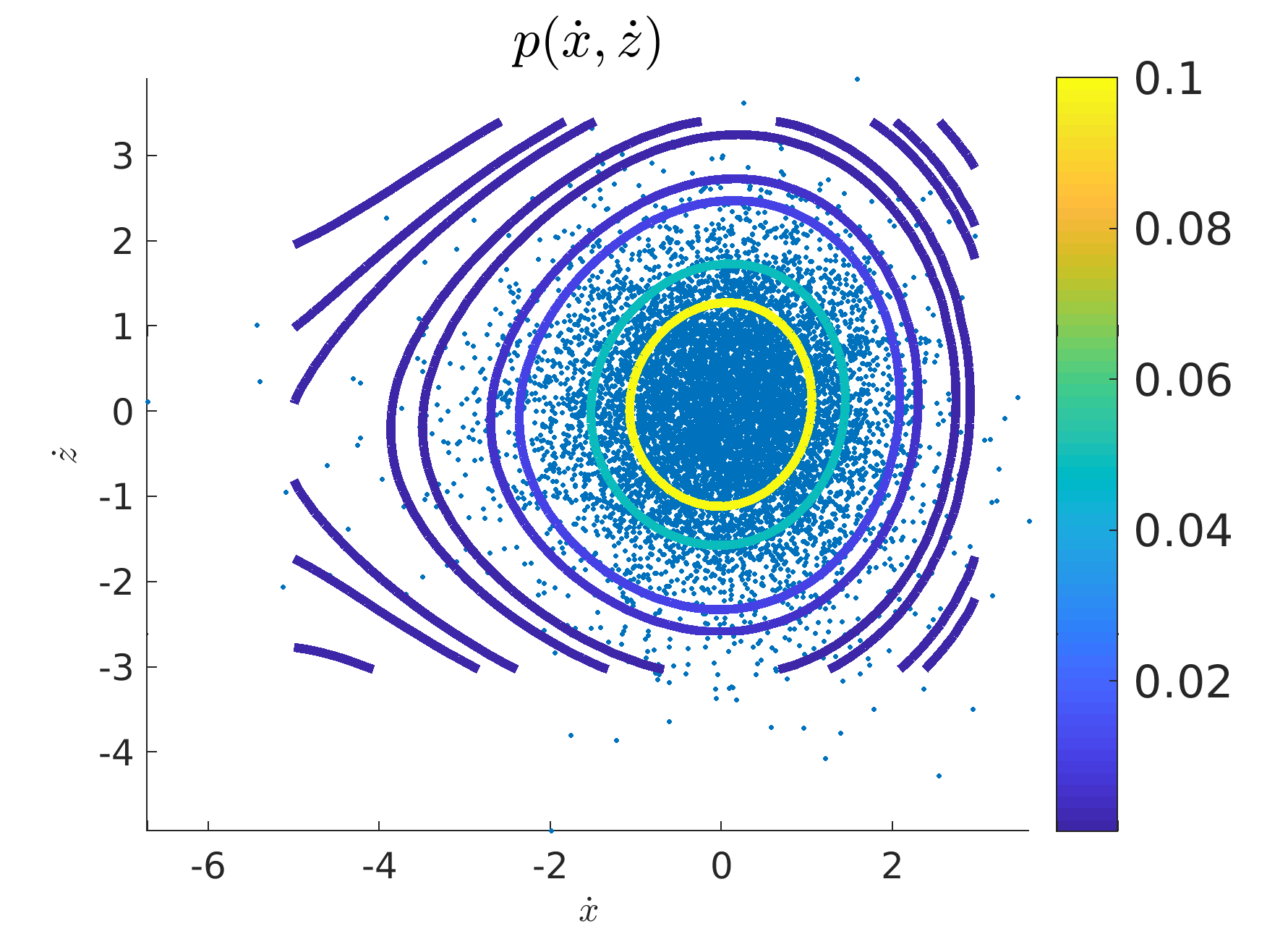}} \\ 
{\multirow{-6}{*}{\rotatebox[origin=c]{90}{\scriptsize  Time = $t_f$  }}} 
& \addheight{\includegraphics[width=1.3in]{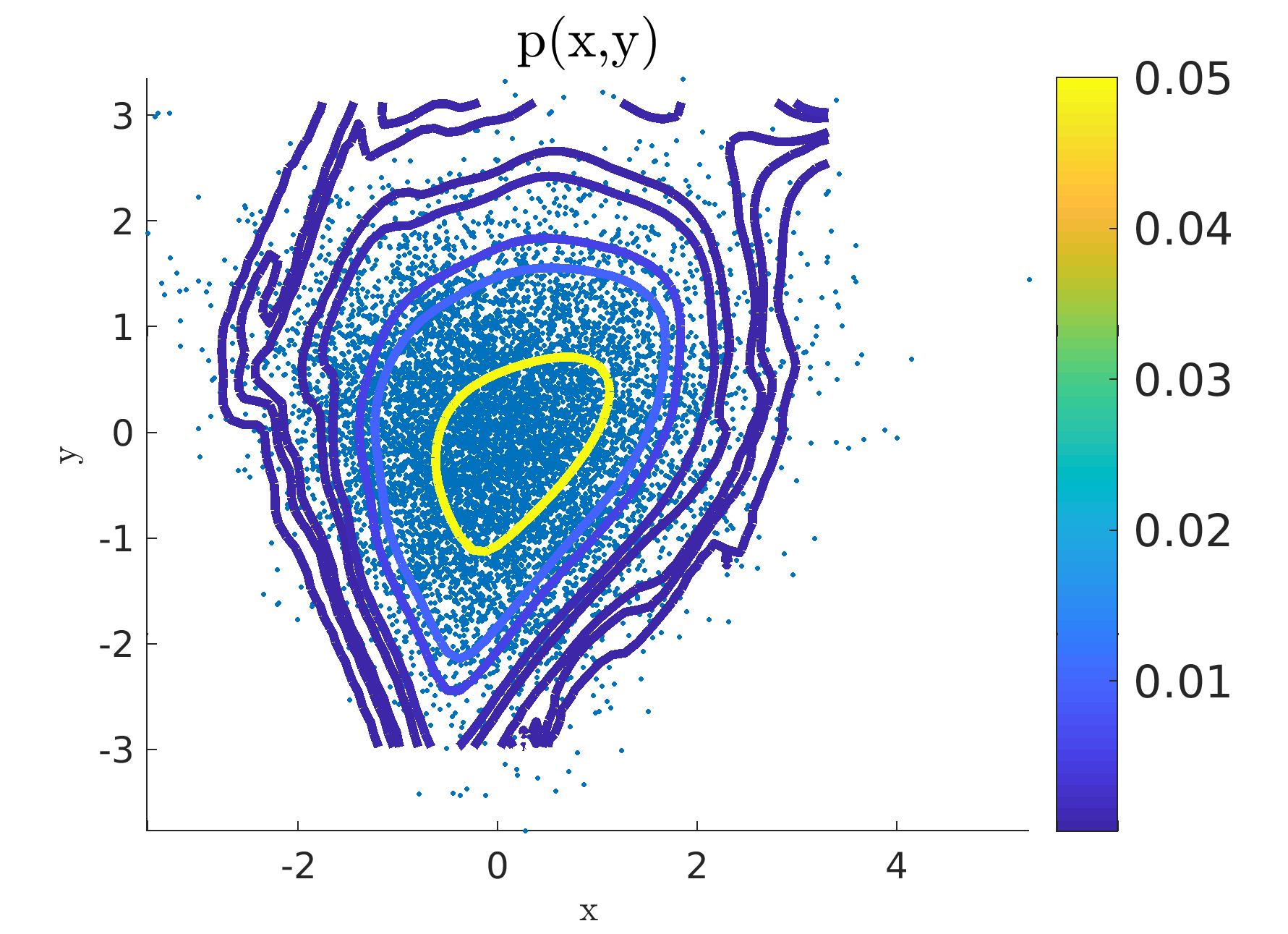}} 
& {\includegraphics[width=1.3in]{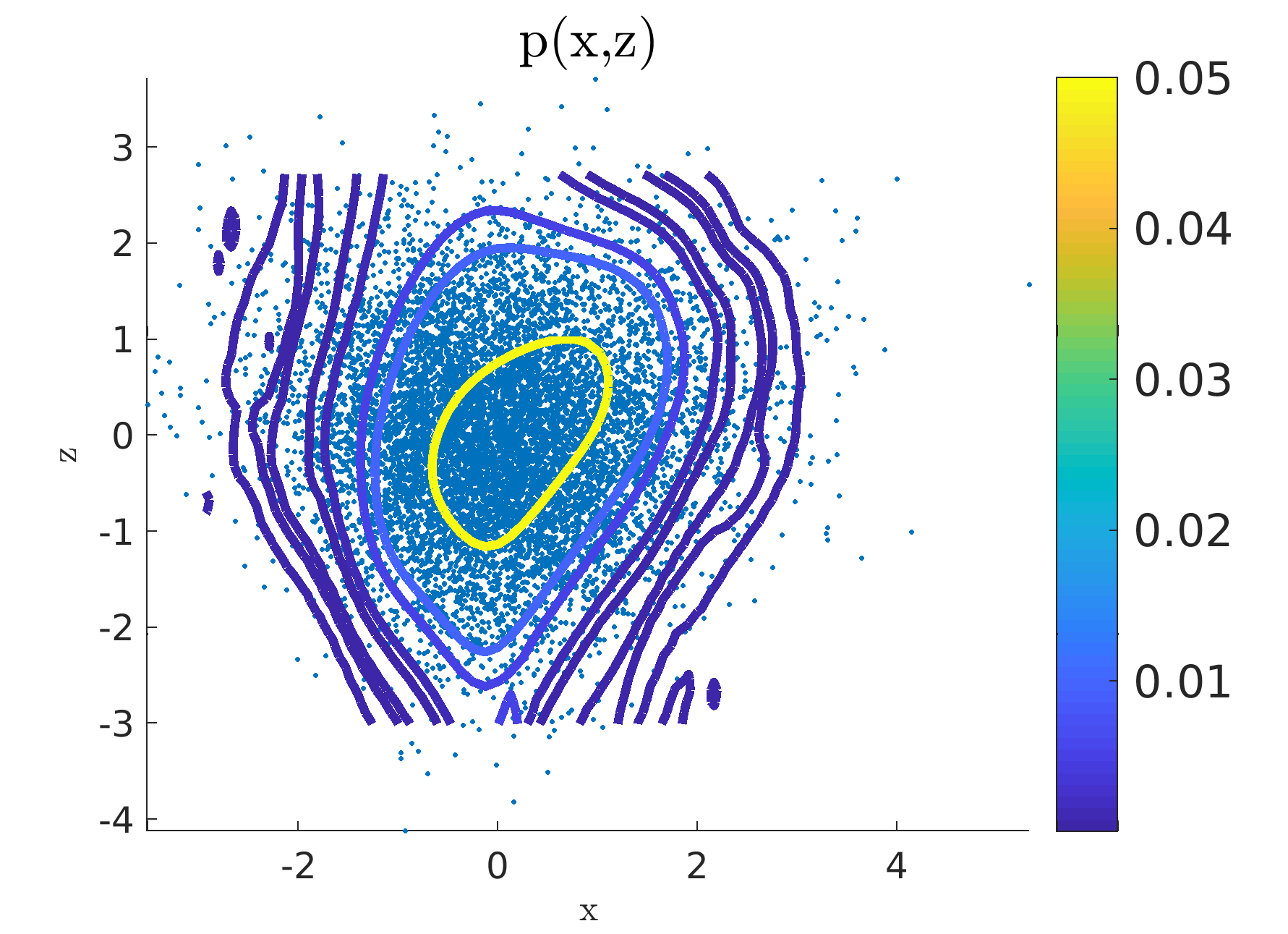}} 
& {\includegraphics[width=1.3in]{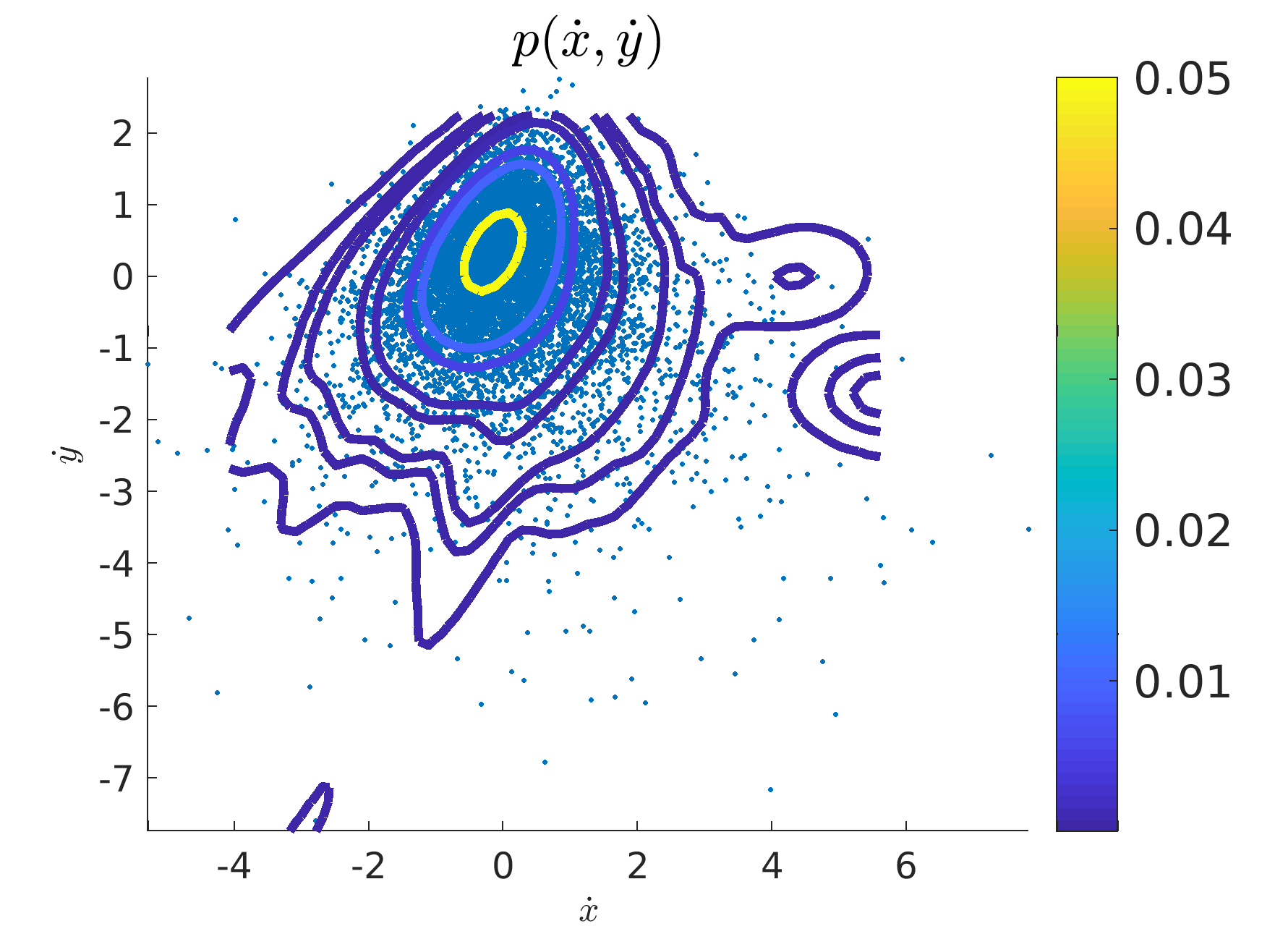}} 
& {\includegraphics[width=1.3in]{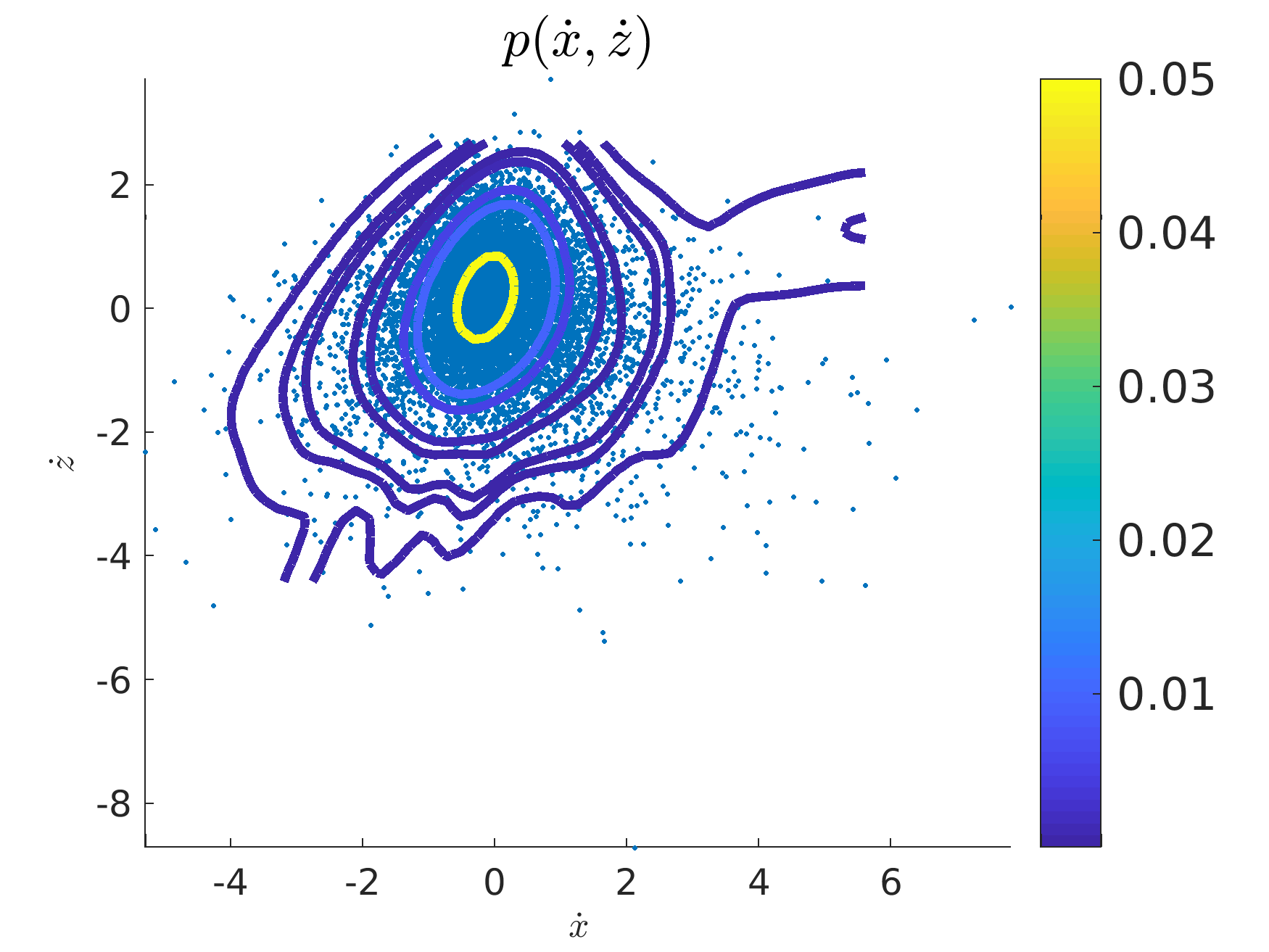}} \\ 
\end{tabular}
\caption{Orbit Transfer Maneuver Case - II - Pdf Contours computed using RS Coefficients (745 CUT points) and superimposed with 50,000 Monte Carlo points} 
\label{fig:TwoBurn_Contour}
\end{figure}

\reffig{fig:beta_ratio} shows the log-pdf ratio, $\beta(H)/\beta(\x)$ of Hamiltonians and monomials in case - II.  As time increases, the log-pdf ratio increases with time, signifying the increase in the dominance of Hamiltonian basis functions in propagating the uncertainty. Finally, \reffig{fig:TwoBurn_Contour} depicts the marginalized pdf contour plots for case - II computed using the RS coefficients in the normalized zero mean and identity covariance space at time $t_f$ and $t_f/2$, where $t_f$ represents the time before the second maneuver. The marginalized pdf contours for $p(x,y)$, $p(x,z)$, $p(\dot x,\dot y)$, and $p(\dot x,\dot z)$ are displayed in columns one through four. As ground truth, the contour plots in \reffig{fig:TwoBurn_Contour} are superimposed with $50,000$ MC points depicted as blue dots. The color bars in the contour plots represent the pdf value at various contour levels. Due to the highly non-Gaussian nature of the pdf, it can be noticed that the pdf contours are coming out to be non-smooth due to the marginalization errors. However, when compared to MC distribution, the proposed method accurately captures the overall shape.

%% file: src/04_conclusion.tex
This research proposed a numerical framework to compute the uncertainty propagation for a non-linear dynamical system. In this study, the uncertainty is quantified using probability density functions (pdfs), which represent the system's uncertainty. The pdf is propagated using the Fokker-Planck-Kolmogorov Equation (FPKE), a linear partial differential equation. Furthermore, a sparse-based collocation method is utilized to approximate the log-pdf from an extensive over-complete dictionary of basis functions. The main contribution of this research is the inclusion of the Hamiltonian function in the basis function dictionary along with the monomials. The monomials are used to approximate the transient behavior of the log-pdf while the Hamiltonians govern the stationary pdf. 

Initial condition uncertainty propagation through nonlinear Duffing oscillator is considered to benchmark the developed approach. The computed solution exactly reproduces the analytical solution for the stationary FPKE. Additionally, the problem of orbit transfer maneuver for satellite motion from LEO to a specific location in space is considered for which no closed-form solution is available. By including the Hamiltonian in the dictionary of basis functions in this problem, it can be conclude that the pdf is computed in a much more accurate and efficient manner. 

% to propagate the uncertainty in the pdf through the dynamical system. The Conjugate Unscented Transformation (CUT) method is used to achieve a minimal representation of collocation points that represent the given domain. Furthermore, a $l_1$-norm optimization is proposed to create a minimal representation for the required coefficients in order to capture the propagated pdf accurately. Combining these two minimal representations allows for more efficient computation of the reachability set. Two representative examples efficiently demonstrate the utility of the proposed approach in computing the reachability analysis.